\documentclass[12pt, oneside]{amsart}
\usepackage[margin=1in]{geometry}

\usepackage{amssymb}

\usepackage{graphicx}

\usepackage{soul}
\usepackage{hyperref}
\usepackage{tikz-cd}

\newtheorem{theorem}{Theorem}[section]
\newtheorem{lemma}[theorem]{Lemma}
\newtheorem{proposition}[theorem]{Proposition}
\newtheorem{corollary}[theorem]{Corollary}

\theoremstyle{definition}
\newtheorem{definition}[theorem]{Definition}
\newtheorem{example}[theorem]{Example}

\theoremstyle{remark}

\numberwithin{equation}{section}

\newcommand{\ed}{\text{ed}}
\newcommand{\GL}{\text{GL}}
\newcommand{\Aff}{\mathrm{Aff}}
\newcommand{\SO}{\text{SO}}
\newcommand{\Orth}{\text{O}}
\newcommand{\R}{\mathbb{R}}
\newcommand{\C}{\mathbb{C}}
\newcommand{\Z}{\mathbb{Z}}
\newcommand{\Aut}{\text{Aut}}
\newcommand{\Diag}{\text{Diag }}
\newcommand{\Stab}{\text{Stab}}
\newcommand{\Perm}{\text{Perm}}
\newcommand{\sgn}{\text{sgn}}
\newcommand{\osym}{S^{\pm}}

\newcommand{\id}{\text{id}}

\newcommand{\lra}{\longrightarrow}
\newcommand{\mdim}{\text{mdim}}
\DeclareMathOperator{\mdimc}{\text{mdim}_{\mathbb{C}}}
\DeclareMathOperator{\mdimr}{\text{mdim}_{\mathbb{R}}}
\newcommand{\re}{\text{Re }}
\newcommand{\im}{\text{Im }}

\newcommand{\ol}{\overline}
\newcommand{\fc}{\mathcal{C}}
\newcommand{\fd}{\mathcal{D}}
\newcommand{\fr}{\mathcal{R}}
\newcommand{\fx}{\mathcal{X}}
\newcommand{\fv}{\mathcal{V}}
\newcommand{\fu}{\mathcal{U}}
\newcommand{\fw}{\mathcal{W}}
\newcommand{\fz}{\mathcal{Z}}

\begin{document}

\title[Minimal Representations of Split Group Extensions]{Minimal Faithful Representations of Split Extensions by Abelian Groups}

\author{Charles Daly}
\address{Max Planck Institute for Mathematics in the Sciences, 04103 Leipzig, Germany}
\curraddr{}
\email{charles.daly@mis.mpg.de}
\thanks{This work was supported by Brown University's SPRINT|UTRA Grant [\emph{Group Structures of the Rubik's Cube and Pyramorphix}]. The first author would like to thank Dani Kaufman for her explanations about the group structure of the Rubik's cube.  He would also like to thank both Philip Eberhart and Myles Miller for their insights into the puzzle, particularly Myles for showing him the `bird configuration' of the Rubik's Pyramorphix.  The second author thanks his algebra teachers for their instruction and guidance over the course of his undergraduate education. He also thanks his parents for their constant encouragement, their celebration of curiosity, and their steadfast commitment to his studies.}

\author{Justin Kingsnorth}
\address{Department of Mathematics, Brown University, Providence, RI 02912, USA}
\curraddr{}
\email{justin\_kingsnorth@brown.edu}
\thanks{}


\date{\today}

\dedicatory{}

\begin{abstract}
Every finite group admits a representation which is `most efficient' in the sense that the representation is faithful and the dimension is minimal.  We call such representations minimal faithful representations.  The image of this representation retains all the information about the abstract group and can come with some additional geometric structure.  For example, a minimal faithful representation could preserve a symplectic form or act on a fixed basis by permutations.  It is the purpose of this paper to address these ideas for finite semi-direct products of groups split by an abelian group whose splitting homomorphism is faithful.  We calculate the minimal faithful dimensions of such representations over both the complex and real numbers, and we express them in terms of group invariants to argue that in some sense these groups are extensions of tori.  In the real case, the dimension of a minimal faithful representation and the structure it preserves encode the possible types of orbifold singularities of the group in question.  We apply our framework to some families of groups to calculate their minimal faithful dimensions and analyze some of their orbifold singularities.
\end{abstract}

\maketitle

\section{Introduction}\label{sec:intro}
In this paper, we propose an analysis of finite group theory influenced by the study of geometric structures.  By Cayley's theorem, for a fixed field $F$, every finite group $G$ admits a faithful representation over an $F$-vector space of some finite dimension.  There exists a minimal such dimension for which $G$ admits an embedding.  We call such a representation a \emph{minimal faithful representation} and denote the minimal dimension by $\mdim_{F} \, G$. In this paper, we will study minimal faithful representations of finite groups $G = AH$, where $A$ is abelian and $H$ is a complementary subgroup whose action on $A$ is faithful. We will largely be interested in the cases when $F = \C$ or $F = \R$. We show that groups of the form $G = AH$ with faithful splitting homomorphisms admit minimal faithful representations into the generalized permutation group $\C^{m} \rtimes S_{m}$, and are therefore in some sense extensions of tori.  We propose considering this question in greater generality to better understand the links between abstract group theoretic properties of groups $G$, their minimal faithful representations, and the geometric structures which they preserve, or fail to. \\
\\
 The study of minimal dimension extends several questions of historical interest in group theory. Jordan showed in 1877 that if a group $G$ admits a faithful representation of degree $n$, it contains an abelian normal subgroup $A$ such that $[G:A] \le j(n)$ for some function $j$ \cite{jordan}. Various estimates for $j(n)$ have since been proposed, with optimal bounds given in \cite{collins}. There is also a rich body of work on faithful irreducible representations of a finite group, beginning with Burnside in 1911 \cite{burnside}. Gasch\"utz, Weisner, and Shoda each provided criteria for the existence of such representations \cite{gaschutz} \cite{weisner} \cite{shoda}, and Zhmud' similarly determined the number of irreducible components required for a faithful representation \cite{zhmud}. Minimal dimension additionally builds on the study of minimal permutation degree, which began with Johnson in 1971 \cite{johnson}. The permutation degree is known to add under products for a large family of groups \cite{wright} and has been determined for all finite simple groups \cite{mazurovsporadic} \cite{mazurovclassical} \cite{vasilyev}.\\
\\
The most general bounds on minimal dimension to date are due to Moret\'o, who showed that $\mdim_\C\, G \le \sqrt{|G|}$ for all but a small family of 2-groups \cite{moreto}. Most other recent results are driven by the concept of essential dimension $\ed_F \, G$, which originated with Buhler and Reichstein in the context of algebraic geometry \cite{buhler} and has since been extended to a number of other areas \cite{merkurjev}. The essential dimension is of interest over arbitrary fields $F$, and it is always bounded above by $\mdim_{F} \, G$ \cite{berhuy}. This has motivated a number of results on the minimal dimensions of $p$-groups \cite{cernele} \cite{bardestani}, where the minimal and essential dimensions are known to coincide \cite{karpenko}.\\
\\
This paper is split into six sections.  Section \ref{sec:intro} is introductory.  Section \ref{sec:preliminaries} is dedicated to preliminaries regarding semi-direct products and representations.  Section \ref{sec:rag} studies both complex and real representations of finite abelian groups and relates the dimension of minimal faithful representations to their invariant factor decomposition.  Section \ref{sec:rseag} contains the main results of this paper which prove that groups of the form $G = AH$ with faithful splitting homomorphism admit minimal faithful representations into the generalized permutation group, and, relates the minimal faithful dimension of $G$ to certain generating sets of the dual group of $A$.  Section \ref{sec:rrseag} provides analogous results to those of Section \ref{sec:rseag}, but in the context of real representations.  Section \ref{sec:applications} applies our theory to the Rubik's Cube group $R_{3}$, examines Moret\'{o}'s bound in our context, and shows how the relation between the minimal complex and real representations remains mysterious. 


\section{Preliminaries}\label{sec:preliminaries}
It is the purpose of this section to provide the reader with definitions and examples of the objects which we work with throughout this paper.  Here we recall the basics of semi-direct products and their representations.  We also provide the group structure of the $3\times 3$ Rubik's cube group $R_{3}$ which we investigate in Section \ref{sec:applications}.  The remainder of this section is spent defining \emph{real} and \emph{almost-complex} representations.  These are likely familiar objects to readers, but rephrased in the spirit of complex representations with additional structure. 


\subsection{Representations of Semi-Direct Products}\label{ssec:bt}
Here we introduce some notation and terminology for semi-direct products used repeatedly through this paper.  Given two groups $N$ and $H$ and a homomorphism $\phi: H \lra \Aut(N)$ one can form the semi-direct product $G = N\rtimes_{\phi} H$.  We call the homomorphism $\phi$ the \emph{splitting homomorphism}.  The group $G$ fits into the short exact sequence $N \xrightarrow{i} G \xrightarrow{q} H$ where $i(n) = (n,1)$ and $q(n,h) = h$.  This short exact sequence is \emph{split} in the sense that $q$ admits a homomorphic section $\sigma: H \lra G$ given by $\sigma(h) := (1,h)$ satisfying $\sigma\circ q = \id_{H}$.  Abusively, we denote a semi-direct product by its constituent pieces $G = NH$ and suppress $\phi$ when clear.  We now provide two examples which will appear again later in this work.

\begin{example}\label{ex:orientedperm}
Let $S_{m}$ act on $\Z_{2}^{m}$ by the standard permutation action and $\phi: S_{m} \lra \Aut(\Z_{2}^{m})$ be the homomorphism associated to this action.  We form the group $\osym_{m} := \Z_{2}^{m} \rtimes_{\phi} S_{m}$ which we call the \emph{oriented symmetric group}.  This is in part because $\osym_{m}$ is isomorphic to the subgroup of permutations of $\{1^{+},1^{-},\hdots, m^{+},m^{-}\}$ which send pairs of the same number to pairs of the same number, but possibly with positive and negative signs reserved.  The subgroup $\Z_{2}^{m} < \osym_{m}$ are the permutations which swap $i^{+}$ and $i^{-}$ for all $1 \leq i \leq m$, whereas $S_{m} < \osym_{m}$ is the subgroup of permutations which send pairs $\{i^{+},i^{-}\}$ to $\{j^{+},j^{-}\}$ taking positives to positives and negatives to negatives.  In this way, if $S_{m}$ acts on a set $X$, then $\osym_{m}$ inherits an action on $X^{\pm}$, the set of positive and negative labelings of $X$.  We call elements of $\osym_{m}$ \emph{oriented permutations}.  
\end{example}  
The next example was the original motivation for this paper.  The $3\times 3$ Rubik's cube group $R_{3}$ is a celebrated mathematical object, and the original goal of this work was to calculate $\mdimc R_{3}$ and $\mdimr R_{3}$.  As it so happens, $R_{3}$ is of the form $AH$ where $A$ is an abelian group and $H$ acts on $A$ faithfully.  Before stating its group structure, we introduce some notation.  For each $m, k \geq 1$, define $\Z_{k,0}^{m}$ to be the kernel of the map $\Z_{k}^{m} \lra \Z_{k}$ which assigns a $m$-tuple of elements in $\Z_{k}$ its sum in $\Z_{k}$.    
\begin{example}\label{ex:rubiks}
Let $R_{3}$ denote the Rubik's Cube group, which consists of the symmetries of the puzzle in Figure \ref{fig:3x3_setup}.  We call the pieces of the puzzle \emph{cubelets} and label the cubelets of the $3 \times 3$ cube by the numbers and letters below.  Cubelets with numbers and letters are called \emph{corner cubelets} and \emph{edge cubelets} respectively.  
\begin{figure}[h]
\centering
\begin{minipage}{0.45\textwidth}
\centering
\includegraphics[width = 4.5cm]{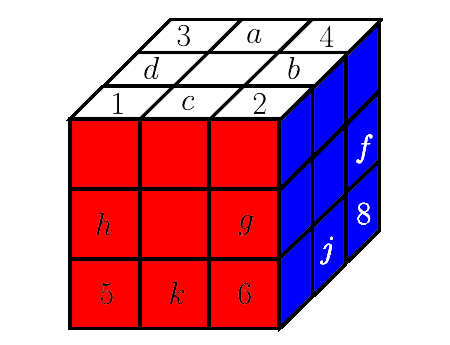}
\end{minipage}
\begin{minipage}{0.45\textwidth}
\centering
\includegraphics[width = 4.5cm]{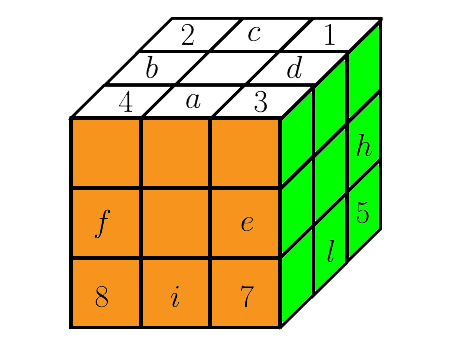}
\end{minipage}
\caption{The $3 \times 3 \times 3$ Rubik's cube, with labels for the corner and edge positions.}
\label{fig:3x3_setup}
\end{figure}
To study the group structure of $R_{3}$, we assign a \emph{local orientation} to each corner position $1 \leq i \leq 8$. This is an element $s_i \in \Z_3$ defined by an orthonormal right-handed basis $\hat{x}_i, \hat{y}_i, \hat{z}_i$ attached to the cubelet at that position.  We color the basis vectors $\hat{x}_i, \hat{y}_i, \hat{z}_i$ yellow, orange, and green respectively.  These bases travel and rotate with the cubelets as they are acted on by $R_{3}$. The number $s_i \in \Z_{3}$ is the number of counterclockwise rotations by $2\pi / 3$ taking the basis at position $i$ in the solved state to the basis at $i$ in the current state.\\
\\
Similarly, we assign each edge position $x \in \{a,b,\hdots, l\}$ an element $t_x \in \Z_2$ defined by an orthonormal basis $\hat{v}_x$, $\hat{w}_x$ attached to the cubelet at that position. We color $\hat{v}_{x}$ pink and leave $\hat{w}_{x}$ implicit.  The number $t_x \in \Z_{2}$ is the number of rotations by $\pi$ taking the current basis at $x$ to the basis in the solved state.\\
\\
Let $A$ be the subgroup of $R_{3}$ which keeps the cubelets in place.  It is clear $A$ is abelian, and to each such move, we assign a pair of tuples $T:=(t_{a},t_{b},\hdots, t_{l}) \in \Z_{2}^{12}$ and $S:=(s_{1},\hdots, s_{8}) \in \Z_{3}^{8}$.  This element may be identified with the number of $\pi$-rotations of edge cubelets and $2\pi/3$-counterclockwise rotations of corner cubelets to be performed to get from the solved state to its current state.  A remarkable fact about the puzzle is that the sums $t_{a} + \hdots + t_{l}$ and $s_{1} + \hdots + s_{8}$ are invariant under $R_{3}$, and by comparison to the solved state, must always be zero.  We may thus identify tuples $(T,S) \in \Z_{2,0}^{12}\oplus \Z_{3,0}^{8}$.  An example of a move $a \in A$ is illustrated in Figure \ref{fig:3x3_localorients}.  
\begin{figure}[h]
\centering
\includegraphics[width = 10cm]{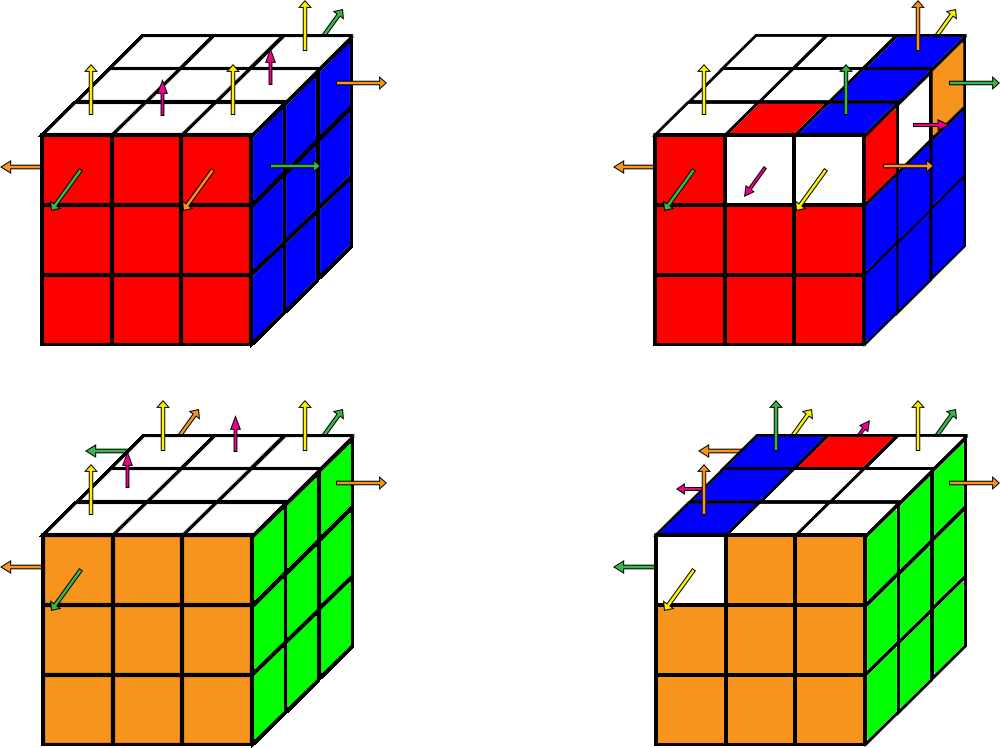}
\caption{At left, bases for some local orientations along the top face.  At right, the change in local orientation under a an element in $\Z_{2,0}^{12}\oplus \Z_{3,0}^{8}$. The changes in local orientations are $(a,b,c,d) = ([0], [1], [1], [0])$, $(e,f,g,h) = ([0], [0], [0], [0])$, $(i,j,k,l) = ([0], [0], [0], [0])$; $(1, 2, 3, 4) = ([0], [1], [0], [2])$, $(5, 6, 7, 8) = ([0], [0], [0], [0])$.
}
\label{fig:3x3_localorients}
\end{figure}
Next consider the subgroup of moves $H$ which preserve the local orientations.  That is, moves which permute the cubelets, but the yellow, orange, green, and pink arrows must all be sent to one another.  An example of an $h \in H$ would be a twist of the top face of the cube relative to the local orientations as in Figure \ref{fig:3x3_twist}.
\begin{figure}[h]
\centering
\begin{minipage}{0.45\textwidth}
\centering
\includegraphics[width = 6cm]{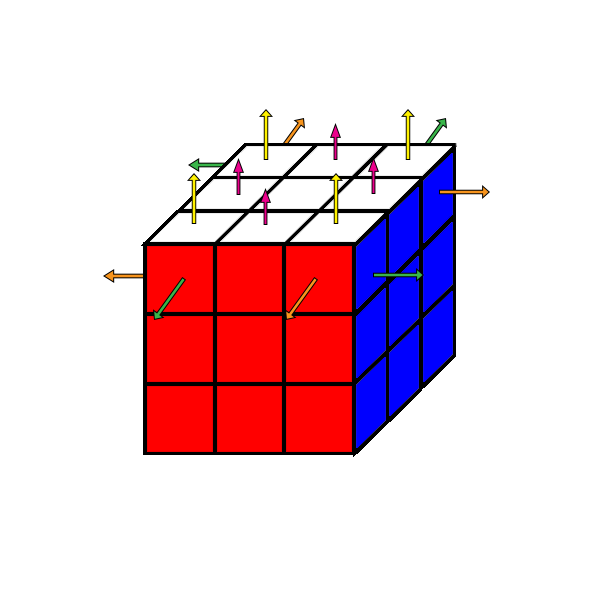}
\end{minipage}
\begin{minipage}{0.45\textwidth}
\centering
\includegraphics[width = 6cm]{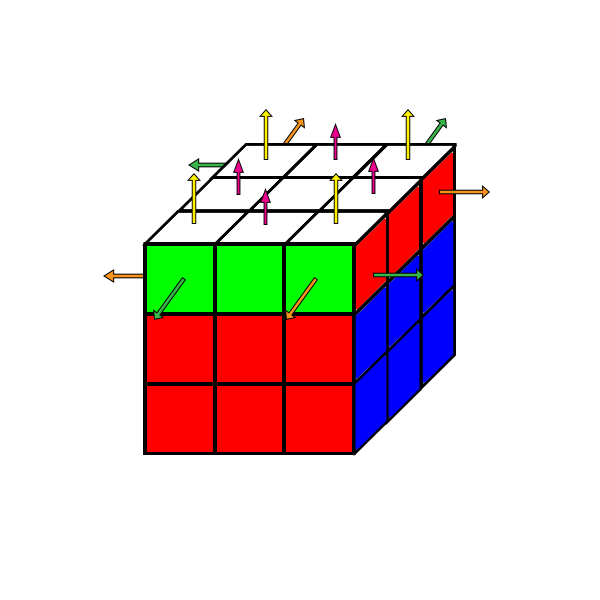}
\end{minipage}
\vspace{-1cm}
\caption{A twist of the cube which preserves the local orientations prescribed above.}
\label{fig:3x3_twist}
\end{figure}
Such a move is therefore determined by a permutation of the cubelets which we identify with $(\alpha,\beta) \in S_{12}\times S_{8}$.  It can be shown that the product $\sgn(\alpha)\sgn(\beta)$ is preserved by all of $R_{3}$, and again by comparison to the identity element, must always be trivial.  
\\
\\
It is clear by definition that $A\cap H = 1$ because if a move fixes each cubelet and does not rotate or flip any of the cubelets, it must be trivial.  Showing that $R_{3}\cong AH$ where $H$ acts on $A$ by permutations is a fact we will take for granted, but use later in Section \ref{sec:applications}.
\end{example}
Recall a finite dimensional representation of a group $G$ over a base field $F$ is a homomorphism from $G$ into the general linear group of $V$, $\GL(V)$, where $V$ is a finite dimensional vector space over $F$.  When the context is clear, we will simply refer to $V$ as the representation.  For our purposes, $F$ will always be $\C$; however, we will reinterpret the case where $F =\R$ as a special case where $F = \C$ and $G$ acts on $V$ by preserving additional structure.  We are interested in representations of semi-direct products $G = NH$, so it is natural to ask if a representation of $G$ is determined by representations on its pieces $N$ and $H$.  This is the content of our first lemma.
\begin{lemma}\label{lem:semi_rep}
Representations $V$ of the semi-direct product $G = NH$ are in one-to-one correspondence with pairs of representations $\rho: N \lra \GL(V)$ and $\sigma: H \lra N_{\GL(V)}\rho(N)$ satisfying the compatibility condition below for all $n \in N$ and $h \in H$.  
\begin{equation}\label{eq:comp-rel}
\left(\rho\circ \phi_{h}\right)(n) = \sigma(h) \rho(n) \sigma(h)^{-1}
\end{equation}
\end{lemma}
\begin{proof}
Let $\rho: G \lra \GL(V)$ be a representation of $G$.  If we abusively denote $\rho$ on $N$ by $\rho$, and $\rho$ on $H$ by $\sigma$, routine calculation shows these representations satisfy the equality in Equation \ref{eq:comp-rel} for all choices of $n$ and $h$.  Conversely, if we begin with such a pair of representations, we can extend the representation $\rho$ on $N$ to all of $G = NH$ by defining $\rho(g) = \rho(nh) := \rho(n)\sigma(h)$.  This is well defined as the decomposition of $g = nh$ is unique.  That it is a homomorphism is due to the calculation below.
\begin{align*}
\rho(gg') &:=  \rho\left(n\phi_{h}(n')\right)\sigma(hh') = \rho(n)\left(\sigma(h)\rho(n')\sigma(h)^{-1}\right)\sigma(h)\sigma(h')\\
&= \left(\rho(n)\sigma(h)\right)\left(\rho(n')\sigma(h')\right) = \rho(g)\rho(g')  
\end{align*}
To see this correspondence is a bijection, note the value of the representation on $g = nh$ is entirely determined by its values on $n$ and $h$ respectively. 
\end{proof}
If we assume the splitting homomorphism is faithful, this imposes some restrictions on the representation, and can be leveraged to guarantee faithfulness of the representation as seen in the lemma below.
\begin{lemma}\label{lem:avoids_central}
Let $G = NH$ be a semi-direct product with faithful splitting homomorphism $\phi$.  If $\rho: N \lra \GL(V)$ and $\sigma: H \lra N_{\GL(V)}\rho(N)$ satisfy the compatibility condition of Lemma \ref{lem:semi_rep} and $\rho$ is faithful, then $\sigma$ intersects the centralizer of $\rho(N)$ trivially and is therefore faithful. If $\rho$ and $\sigma$ satisfy the compatibility condition and the image of $\sigma$ intersects the image of $\rho$ trivially, then the pair $(\rho, \sigma)$ determines a faithful representation of $G$ if and only if both $\rho$ and $\sigma$ are faithful.
\end{lemma}
\begin{proof}
If $\rho$ and $\sigma$ satisfy the compatibility condition, they determine a representation of $G$ which we also denote by $\rho$.  We claim the restriction of $\rho$ to $H$, which is $\sigma$ by definition, avoids the centralizer of $\rho(N)$.  For if there is an $h \in H$, where $\rho(h) = \sigma(h)$ is in the centralizer of $\rho(N)$, then for every $n \in N$ we have the equality below. 
\begin{equation*}
\rho(n) = \sigma(h) \rho(n) \sigma(h)^{-1}  = \rho(hnh^{-1}) = \rho\left( \phi_{h}(n)  \right)
\end{equation*}
By faithfulness of $\rho$, $\phi_{h}(n) = n$ for all $n \in N$, and by the faithfulness of the splitting homomorphism, $h = 1$ proving our claim that $\sigma$ intersects the centralizer of $\rho$ trivially.\\
\\
Assume in addition the image of $\sigma: H \lra N_{\GL(V)}(\rho(N))$ intersects the image of $\rho: N \lra \GL(V)$ trivially.  We prove the representation determined by the pair $(\rho,\sigma)$ is faithful if and only if both $\rho$ and $\sigma$ are faithful on $N$ and $H$ respectively.  One direction is clear.  Let us assume both $\rho$ and $\sigma$ are faithful.  If there is a $g = nh \in G$ so that $\rho(g) := \rho(n)\sigma(h) = 1$, then by the fact that $\sigma$ intersects the image of $\rho$ trivially, we have $\sigma(h) = \rho(n) = 1$.  Thus $ h = n = 1$ so $g = 1$.  
\end{proof}

\subsection{Real and Almost-Complex Representations}\label{ssec:racr}

We are motivated to study both complex and real representations of $G$, the latter of which find motivation in the study of geometric structures, otherwise known as $(G,X)$-structures.  A bit more precisely, one endows a manifold $M$ with a local structure of a model space $X$ via charts in such a way that the coordinate changes are given by restrictions of symmetries in $G$ of $X$.  This process gives rise to a \emph{holonomy representation} of the fundamental group of $M$ into the group $G$.  For example, in the study of locally homogenous Riemannian surfaces, the corresponding geometries are spherical, Euclidean, and hyperbolic with model spaces given by the sphere, the Euclidean plane, and the upper half-plane respectively.  The groups of orientation preserving isometries of these geometries are real algebraic Lie groups isomorphic to $\text{SO}(3,\R), \R^{2}\rtimes\text{SO}(2,\R)$ and $\text{PSL}(2,\R)$ respectively.  \\
\\
While finite groups are not typically not the fundamental groups of manifolds of non-negative curvature, they can appear as subgroups of stabilizers of $G$ acting on $X$.  The quotients of the model space by these subgroups model $(G,X)$-\emph{orbifold singularities}.  A simple example is given by the cyclic subgroup of rotations by $2\pi/n$ radians about a point in $\R^{2}$.  The corresponding orbifold singularity about the fixed point is a \emph{cone}.  An very accessibly written treatment of orbifolds can be found in Caramello Jr.'s work \cite{Caramello2019Introduction}.  
\\
\\
If $H \subset G$ is a finite subgroup of a stabilizer of $x \in X$, then we may choose charts so that $H$ is locally acting linearly near $x$ and thus $H$ may be thought of a subgroup of $\Orth(n,\R)$ where $n$ is the dimension of $X$.  Thus by studying the minimal real representations of $H$, we are implicitly studying the minimal real dimension of the corresponding orbifold singularity.  By studying the structure of minimal real representations of a finite group and what type of geometric structures they preserve, or fail to, we gain insight into orbifold singularities of the group in question.  Explicit examples of these dimensions and structures are provided in Section \ref{sec:applications}.
\\
\\
As we wish to study types of representations, we introduce a framework for working with real representations as representations over the complex numbers that preserve additional structure. In what follows all vector spaces are complex, and all maps are assumed to be complex linear unless otherwise explicitly stated.  Much of what follows in this section is inspired by Poonen \cite{realreps}.  Before proceeding, we briefly introduce some notation that will be used throughout this paper.  \\
\\
Given a complex representation $V$ of $G$ we denote by $\ol{V}$ the conjugate representation of $G$.  The vectors of $\ol{V}$ are the same as those of $V$, but the scalar multiplication on $\ol{V}$ is defined to be $c\cdot v := \ol{c}v$ for all $c \in \C$ and $v \in \ol{V}$.  We say an additive map $J: V \lra V$ is \emph{conjugate-linear} if for all $c \in \C$ and $v \in V$, $J(cv) = \overline{c}J(v)$, or equivalently, $J: V \lra \ol{V}$ is linear.  An \emph{involution} is a linear map whose square is $\id$ and an \emph{anti-involution} is a map whose square is $-\id$.  For a map $g: V \lra V$, we denote by $V^{g}$ the subspace of $g$-invariants of $V$, and $E_{\lambda}^{g}$ to be the $\lambda$-eigenspace of $g$.   
\\ \\
With these conventions, we define real representations as follows.

\begin{definition}\label{def:real_rep}
A \emph{real representation} $(W, J)$ of $G$ is a complex representation $W$ of $G$ where $G$ preserves a conjugate-linear involution $J$ in the sense that $JgJ^{-1} = g$ for all $g \in G$.  Alternatively a real representation is a homomorphism from $G$ into the group $\GL(W,J)$ of all automorphisms of $W$ which preserve $J$.  We call $J$ a \emph{real structure}.  
\end{definition}  

This definition is equivalent to the traditional formulation. Given a real representation $(W, J)$ as defined above, we may obtain a representation over the real numbers by restricting scalars to the real vector subspace of $J$-invariants, $W^J$. Conversely, given a representation $U$ over the real numbers, one may obtain a complex representation $W := \C \otimes_\R U$ where $g \in G$ acts by $g(z \otimes u) := z \otimes gu$. The involution $J(z \otimes u) = \ol{z} \otimes u$ then supplies a real representation in our sense. 
\\ \\
Following Definition \ref{def:real_rep}, a \emph{real} $G$-module homomorphism between two real representations $(W,J)$ and $(W',J'$) is a complex $G$-module homomorphism $\phi: (W,J) \lra (W',J')$ which satisfies the property that $\phi\circ J = J' \circ \phi$. The notions of isomorphisms, embeddings, and quotients of real representations follow in a similar way. We will often drop the adjective \emph{real} when the context is clear, and we will typically reserve $W$ for real representations and $V$ for complex representations.
\\ \\
We also introduce almost-complex representations, which admit a variant of complex scalar multiplication.
\begin{definition}\label{def:almostcomplex_rep}
An \emph{almost-complex} representation is a triple $(W,J,I)$ where $(W,J)$ is a real representation, and $I$ is a $G$-invariant anti-involution $I: W \lra W$ that is compatible with the real structure $J$ in the sense that $JIJ^{-1} = I$.  We say a real representation $(W,J)$ admits an almost-complex structure if there is such an $I$ for which $(W,J,I)$ is an almost-complex representation.  
\end{definition}

One may construct an almost-complex representation from any complex representation $V$ of $G$. Given such a representation, let $W = V \oplus \ol{V}$, and let $J(v_1, v_2) = (v_2, v_1)$. Define an anti-involution $I: W \lra W$ by $I(v_1, v_2) := (i v_1, -i \cdot v_2)$, recalling that $i \cdot v_2 = -i v_2$ by our conventions for multiplication in $\ol{V}$. One can show that the triple $(W, J, I)$ is in line with Definition \ref{def:almostcomplex_rep}. 
\\ \\
In fact, all almost-complex representations arise in this way.
\begin{proposition}\label{lem:almostcomplexeq}
A real representation $(W,J)$ admits an almost-complex structure if and only if there exists a complex representation $V$ such that $(W ,J)$ is isomorphic to the real representation $(V\oplus \overline{V},J')$ where $J'(v_{1},v_{2}) = (v_{2},v_{1})$.  
\end{proposition}
\begin{proof}
For the reverse direction, follow the construction above. For the forward direction, let $V$ and $\ol{V}$ be the eigenspaces $E_i^{I}$ and $E_{-i}^I$ of the $G$-invariant anti-involution $I: W \lra W$.
\end{proof}
Not all real representations admit almost-complex structures. For example, the usual representation of $\SO(2, \R)$ is almost-complex, while the usual representation of $\Orth(2, \R)$ is not. 
\\ \\
To close the section, we formally define subrepresentations and irreducible representations over the real numbers, and we state analogues to Schur's lemma and Maschke's theorem in this context. We omit the proofs, which follow with small modifications from the standard arguments.

\begin{definition}\label{def:submodule}
A \emph{real subrepresentation} $(W',J|_{W'})$ of $(W,J)$ is a $J$-invariant complex subrepresentation $W' \leq W$.  We say a real representation of $(W,J)$ is \emph{irreducible} if $(W,J)$ admits no real subrepresentations.
\end{definition}
\begin{lemma}[Schur's lemma]\label{lem:schur}
Let $(W,J)$ be an irreducible real representation and let $\phi: (W,J) \lra (W,J)$ be a real $G$-module homomorphism.  Either $\phi$ acts by scalar multiplication by a real number, or, there exists a complex number $\lambda$ with non-zero imaginary part and an almost-complex structure $I: W \lra W$ such that $\phi$ acts by $\emph{\re}(\lambda)\id + \emph{\im}(\lambda)I$.  
\end{lemma}
\begin{lemma}[Maschke's theorem]\label{lem:maschke}
Let $G$ be a finite group and $(W,J)$ be a real representation of $G$.  Then $W$ splits into the sum of irreducible real representations, $W = W_{1} \oplus \hdots \oplus W_{m}$.
\end{lemma}

\section{Representations of Abelian Groups}\label{sec:rag}
It is the purpose of this section to provide some results on faithful representations, both complex and real, of finite abelian groups $A$.  The main results of this section are Lemma \ref{lem:normalform} wherein we provide a normal form for representations of $A$, and Theorem \ref{thm:min_abl} where we construct minimal faithful representations of $A$.  The latter result shows that, in some sense, all finite abelian groups are toric as they admit minimal faithful representations into tori.  
\\
\\ 
We first introduce some notation.  We denote by $S^{1}$ the group of complex numbers $S^{1} = \{z \in \C \, | \, |z| = 1\}$ under multiplication.  Each character $\chi$ of $A$ is naturally considered as an element of the dual group of $A$, the group of homomorphisms from $A$ to $S^{1}$ which we denote by $\hat{A}$.  We say a character $\eta: A \lra S^{1}$ is \emph{real} if $\eta$ takes all of its values in $\{\pm 1\}$, and say a character is \emph{complex} otherwise.  This classification partitions the group $\hat{A}$ and a character $\eta$ is readily seen to be real if and only if $\eta^{2} = 1$.  We will typically reserve the symbol $\eta$ for real characters, and let $\chi$ be either type depending on the context.  For any given $k \geq 1$, we identify $\Z_{k} \leq S^{1}$ by the embedding taking the generator $1 \in \Z_{k}$ to $\exp(2\pi i/ k) \in S^{1}$.  Finally, the following definition will be used repeatedly throughout the course of this paper to simplify notation.  
\begin{definition}\label{def:actsbychi}
Let $V$ be a representation of a finite abelian group $A$.  We say that $A$ acts by $\chi \in \hat{A}$ on $V$ if $av = \chi(a)v$ for all $v \in V$ and $a \in A$.  Let $(W,J,I)$ be an almost-complex representation of $A$.  We say that $A$ acts by $\chi \in \hat{A}$ on $(W,J,I)$ if $av = \re{\chi(a)}v + \im{\chi(a)}Iv$ for all $v \in W$ and $a \in A$.  
\end{definition}  
Note that $A$ acts by $\chi$ on $(W,J,I)$ if and only if $A$ acts by $\overline{\chi}$ on $(W,J,-I)$.  This ambiguity is an artifact of the lack of a canonical choice to the equation $I^{2} = -\id$.  However, if we fix an almost-complex structure we may unambiguously say $A$ acts by $\chi$ on $(W,J,I)$.  In so doing, $W$ admits a decomposition as $W = E_{i}^{I}\oplus E_{-i}^{I}$ where $A$ acts by $\chi$ on $E_{i}^{I}$ and by $\ol{\chi}$ on $E_{-i}^{I}$ and $J$ permutes the two pieces of these decompositions as guaranteed by Lemma \ref{lem:almostcomplexeq}.  With these definitions established, we prove our first lemma classifying irreducible real representations of abelian groups.  Note this lemma does not assume finiteness of the abelian group $A$.
\begin{lemma}\label{lem:irr_ab_real}
Let $A$ be an abelian group and $(W,J)$ be an irreducible real representation.  Either $W$ is one-dimensional and $A$ acts by some real character $\eta$, or, $W$ is two-dimensional and $(W,J)$ admits an almost-complex structure for which $A$ acts by a complex character $\chi$.  
\end{lemma}
\begin{proof}
Because $A$ is commutative, every $a \in A$ is a real $A$-module homomorphism of $(W,J)$, so by Lemma \ref{lem:schur}, either $A$ acts on $W$ by a real character $\eta$, or, $(W,J)$ admits an almost-complex structure $I$, and $A$ acts by some complex character $\chi$.    \\
\\
It remains to prove the dimension claim.  Let $v \in W$.  If $Jv = zv$ for some $z \in \C$, then the space spanned by $v$ is one-dimensional and $J$-invariant.  Otherwise the space spanned by $\{v,Jv\}$ is two-dimensional and $J$-invariant as $J$ is a conjugate-linear involution.  As $A$ is abelian, any $J$-invariant subspace of $W$ is a real subrepresentation.  By irreducibility, $W$ must be one or two dimensional. 
\end{proof}
We now have the necessary tools to prove the following lemma which gives us a normal form for both complex and real representations of finite abelian groups $A$.
\begin{lemma}\label{lem:normalform}
Let $A$ be a finite abelian group and let $V$ be a representation of $A$.  There exists a collection of distinct characters $\fx = \{\chi_{1},\hdots, \chi_{m}\}$ of $A$ and a decomposition $V = V_{1}\oplus \hdots \oplus V_{m}$ such that each $a \in A$ acts by $\chi_{i}$ on $V_{i}$.
\\
\\
Let $(Z,J)$ be a real representation of $A$.  There exist collections $\fr = \{\eta_{1},\hdots, \eta_{n}\}$, $\fx = \{\chi_{1},\hdots ,\chi_{m}\}$ of distinct real characters and distinct pairwise non-conjugate complex characters and a decomposition $Z = \left(U_{1} \oplus \hdots \oplus U_{n}\right) \oplus \left(W_{1} \oplus \hdots \oplus W_{m}\right)$ as real representations.  Moreover, there is an almost-complex structure $I$, compatible with $J$, which preserves the decomposition of $W$ so that $A$ acts by $\eta_{i}$ on each $U_{i}$ and by $\chi_{i}$ on each $(W_{i}, J, I)$.   
\end{lemma}
\begin{proof}
The first part of the above lemma follows from both Maschke's Theorem for complex representations and the fact that every irreducible representation of an abelian group is one-dimensional.  This provides a decomposition of $V = U_{1}\oplus \hdots \oplus U_{m}$ where each $U_{i}$ is a one-dimensional irreducible representation of $A$.  Each $U_{i}$ determines a character $\chi_{i}$ where $A$ acts on $U_{i}$ by $\chi_{i}$.  We form the collection of \emph{distinct} characters $\{\chi_{1},\hdots, \chi_{m}\}$ and then define each $V_{i}$ to be the isotypic component of $\chi_{i}$.   \\
\\ 
The proof in the second case follows the same arguments using Lemma \ref{lem:maschke} and Lemma \ref{lem:irr_ab_real}.  These lemmas guarantee $(W,J)$ splits into the sum of one-dimensional subrepresentations where $A$ acts by real characters $\eta_{i}$, and, two-dimensional almost-complex subrepresentations where $A$ acts by complex characters $\chi_{i}$ for almost complex structures $I_{i}$.  The $U_{i}$'s in the decomposition are the isotypic components of $\eta_{i}$ and the $W_{i}$'s are the isotypic components of the $\chi_{i}$'s and their conjugates. 
\end{proof}
In light of Lemma \ref{lem:normalform}, we see that a representation of a finite abelian group $A$ determines a collection of distinct characters $\{\chi_{1},\hdots, \chi_{m}\}$.  Faithfulness of the representation is readily seen to be equivalent to requiring $\bigcap_{i=1}^{i=m} \ker \chi_{i} = 1$.  The lemma below shows this is equivalent to the collection of characters being a generating set.  
\begin{lemma}\label{lem:dual_group_gen}
Let $V$ be a representation (complex or real) of a finite abelian group $A$.  The representation is faithful if and only if the collection of characters $\fx = \{\chi_{1},\hdots, \chi_{m}\}$ generates the whole dual group.  
\end{lemma}
\begin{proof}
Given a subgroup $\hat{B} \leq \hat{A}$, we associate the subgroup $B := \{a \in A \, | \, \chi(a) = 1\text{ for all }\chi \in \hat{B}\}$.  Given a subgroup $B \leq A$, we associate the subgroup $\hat{B} = \{\chi \in \hat{A} \, | \, B \leq \ker \chi\}$.  These correspondences are inverses, and determine a bijection between subgroups of $A$ and subgroups of $\hat{A}$.  The claim follows from this bijection, and the fact that $\bigcap_{\chi \in \hat{A}} \ker \chi = \bigcap_{i=1}^{i=m} \ker \chi_{i}$ if $\{\chi_{1},\hdots, \chi_{m}\}$ if $\fx$ is a generating set.  
\end{proof}
We conclude this section with our first theorem which relates the minimal faithful dimensions of finite abelian groups to their invariant factor decomposition. 
\begin{theorem}\label{thm:min_abl}
Let $A$ be a finite abelian group.  Then $\mdimc A $ is equal to the size of a minimal generating set of $\hat{A}$ whereas $\mdimr A $ is equal to the size of a conjugation invariant minimal generating set of $\hat{A}$.  If the invariant factor decomposition of $A$ is given by $A = \Z_{2}^{a} \oplus \Z_{d_{a+1}} \oplus \hdots \oplus \Z_{d_{a+b}}$ where $d_{1} = \hdots = d_{a} = 2$ and $d_{a+1} > 2$, and $d_{i} | d_{i+1}$ for all $1 \leq i < a+b$, then $\mdimc A = a+b$ and $\mdimr A = a+2b$. 
\end{theorem}
\begin{proof}
Given a minimal generating set $\{\chi_{1},\hdots, \chi_{m}\}$ of $\hat{A}$ form the vector space spanned by the basis $\{e_{\chi_{1}},\hdots ,e_{\chi_{m}}\}$ and have $A$ act on each $V_{i} := \C e_{\chi_{i}}$ by $\chi_{i}$.  By Lemma \ref{lem:dual_group_gen}, this representation is faithful.  We claim this dimension is minimal.  Given a faithful representation of $A$, we have an associated collection of characters $\{\chi_{1}, \hdots , \chi_{m}\} \subset \hat{A}$.  By Lemma \ref{lem:dual_group_gen}, this collection generates $\hat{A}$, hence $m$ is at least as large as the size of a minimal generating set, proving the first claim.  The real analogue follows from the fact that every faithful representation determines a conjugation invariant generating set of $\hat{A}$, and a smallest such one provides a real representation of the desired dimension.   \\
\\
To prove the equalities, we let $A$ have the invariant factor decomposition as stated in the theorem.  As $\hat{A}$ and $A$ are isomorphic as groups, they share the same invariant factor decomposition.  Moreover, the size of a smallest generating set of $A$ is equal to the number of invariant factors.  This proves that $\mdimc A  = a+b$.  \\
\\
Let $S$ be a conjugation invariant generating set of $\hat{A}$.  We will show that $S$ has at least $a+2b$ elements.  Note that $S$ separates into real characters and complex characters, say $S = \{\eta_{1},\hdots ,\eta_{r}\} \cup \{\chi_{1},\ol{\chi_{1}},\hdots, \chi_{c},\ol{\chi_{c}}\}$ where the $\eta_{i}$ are real and the $\chi_{i}$ are complex.  Form the new set $S' = \{\eta_{1},\hdots, \eta_{r}\} \cup \{ \chi_{1},\hdots, \chi_{c}\}$.  The set $S'$ is still a generating set of $\hat{A}$ because $\chi\overline{\chi} = 1$ for each complex character $\chi$.  Thus we have that $|S'| = r+c \geq a+b$.  If we show that $c \geq b$, then $|S| \geq a + 2b$ and we will be done.  \\
\\
To this end, let $C$ be the subgroup of $\hat{A}$ generated by $\{ \chi_{1},\hdots, \chi_{c}\}$.  Consider the quotient map $\hat{A} \lra \hat{A}/\hat{A}^{J}$ where $J$ denotes the involution of conjugation and $\hat{A}^{J}$ are the conjugation invariant elements of $\hat{A}$.  By definition, the subgroup generated by $\{\eta_{1},\hdots, \eta_{r}\}$ is contained in the kernel of this map.  Therefore we have an induced map $C \lra \hat{A}/\hat{A}^{J}$ which is a quotient because $S'$ generates all of $\hat{A}$.  It is clear that $\hat{A}/\hat{A}^{J} \cong \Z_{d_{a+1}/2} \oplus \hdots \oplus \Z_{d_{a+b}/2}$ has $b$ invariant factors, therefore $C$ must have at least $b$ invariant factors.  Thus $c \geq b$ so $\mdimr A = a+2b$.  
\end{proof}

\section{Representations of Split Extensions by Abelian Groups}\label{sec:rseag}
Recall that a group $G$ \emph{acts by permutations} on $V$ if there exists a fixed basis $\{v_{1},\hdots, v_{m}\}$ of $V$ which is invariant under the action of $G$.  That is to say, for each $g \in G$, there exists a $\sigma \in S_{m}$ so that $gv_{i} = v_{\sigma(i)}$ for each $v_{i}$.  Thus, the group of permutations on $V$ relative to a fixed basis is isomorphic to $S_{m}$.  We say a group acts by \emph{generalized permutations} if $G$ preserves a fixed basis, up to scales.  That is to say, for each $g \in G$, there is a $\sigma \in S_{m}$ and scales $\lambda_{i} \in \C^{*}$ so that $gv_{i} = \lambda_{\sigma(i)}v_{\sigma(i)}$.  The group of generalized permutations relative to a fixed basis is isomorphic to the group $(\C^{*})^{m} \rtimes S_{m}$ where $S_{m}$ acts on $(\C^{*})^{m}$ by the standard permutation action.  \\
\\
For the rest of this paper, let $G = AH$ be a finite group where $A$ is abelian and the splitting homomorphism $\phi: H \lra \Aut(A)$ is faithful.  It is the goal of this section to show these groups admit minimal faithful representations into the generalized permutation group, concluding with Theorem \ref{thm:cmdim}.  As we saw in the previous section, it was gainful to study representations of abelian groups by their corresponding characters.  We will employ the same techniques in this section, but in addition keep track of the associated $H$-action on the dual group $\hat{A}$ which is defined by $(h \chi) : A \lra S^{1}$ via $(h\chi)(a) := \chi(h^{-1}ah) = \left(\chi \circ \phi_{h^{-1}}\right)(a)$ for all $a \in A$.  As we will see in Corollary \ref{cor:toremb}, the minimal faithful dimension of $G = AH$ is the minimal dimension in which we can represent this action of $H$ by permutations.  Before addressing the representation theory of such $G$, we provide a lemma which will help us in our constructions.   
\begin{lemma}\label{lem:dual_action_faithful}
Let $G = AH$ be as above, and $\fx = \{\chi_{1},\hdots, \chi_{m}\}$ be an $H$-invariant generating set of $\hat{A}$.  If $h\chi_{i} = \chi_{i}$ for all $1 \leq i \leq m$, then $h = 1$.  In particular, this means the homomorphism permuting characters of $\fx$, $H \lra S_{m}$, is faithful.  
\end{lemma}
\begin{proof}
If we found such an $h \in H$, it must stabilize all of $\hat{A}$ because $\fx$ generates $\hat{A}$.  If we fix an $a \in A$, then $(h\chi)(a) = \chi(a)$ for all $h \in H$ and any $\chi \in \hat{A}$.  Thus $\chi\left(h^{-1}ah\right) = \chi(a)$ for all $\chi \in \hat{A}$.  If two elements $a,b \in A$ evaluate to the same quantity on every character in $\hat{A}$, then $a = b$.  Hence $\phi_{h^{-1}}(a) = a$ for all $a \in A$, and by faithfulness of the splitting homomorphism, $h = 1$.   
\end{proof}
With this lemma established, we move onto the representation theory of such groups.  Here we introduce notation that will persist throughout this section.  Let $V = V_{1}\oplus \hdots \oplus V_{m}$ be a vector space where $d_{i} := \dim V_{i}$ and denote the ordered list of (possibly repeated) dimensions by $\mathcal{D} := (d_{1},\hdots, d_{m})$.  Let $\fx := \{\chi_{1},\hdots,\chi_{m}\}$ be an $H$-invariant collection of \emph{distinct} characters of the finite abelian group $A$.  Let $D := \Diag V(\fx) $ be the group of automorphisms of $V$ which act on each $V_{i}$ by $\chi_{i}$.  Let $P := \Perm(\fd) \leq S_{m}$ be the subgroup of all permutations $\sigma$ of $\{1,\hdots, m\}$ where $\sigma(i) = j$ implies $d_{i} = d_{j}$, i.e. the subgroup of all permutations of $S_{m}$ that permute elements $\{1,\hdots, m\}$ with the same corresponding dimensions $d_{i} = d_{j}$.  Let $\GL(V)$ be the entire group of automorphisms of $V$, and finally, let $S:= \Stab(\fv)$ be the subgroup of $\GL(V)$ which preserves the decomposition $\fv$ of $V$ into $V_{1}\oplus \hdots \oplus V_{m}$.  
\begin{lemma}\label{lem:complex_diag_norm}
Adopting the notation in the above paragraph, there exists a surjective homomorphism $q$ from the normalizer of $\Diag V(\fx)$ inside of $\GL(V)$ onto $\Perm(\fd) \leq S_{m}$ whose kernel is $\Stab(\fv)$ which is in the centralizer of $\Diag V(\fx)$.  Moreover, this short exact sequence is split.  Symbolically, the exact sequence in Equation \ref{eq:cdn_split} is split.  
\begin{equation}\label{eq:cdn_split}
\Stab(\fv) \xrightarrow{i} N_{\GL(V)}\left(\Diag V(\fx) \right) \xrightarrow{q} \Perm(\fd)
\end{equation} 
\end{lemma}
\begin{proof}
To make notation less cumbersome, let $S, D, N,$ and $P$ denote the stabilizer, diagonal group, normalizer of the diagonal group, and permutation group appearing in Equation \ref{eq:cdn_split}.  We first construct a homomorphism $q: N \lra P$.  We claim each $n \in N$ must permute the pieces of the decomposition $V = V_{1}\oplus \hdots \oplus V_{m}$.  If so, then we define $q(n) \in P \leq S_{m}$ by the permutation $n$ induces on $\{V_{1},\hdots, V_{m}\}$.\\
\\ 
To see that $N$ permutes the $V_i$, fix some $n \in N$ and $1 \le i \le m$. The space $V_i$ lies within an eigenspace of every $a \in A$, as $a$ acts on $V_{i}$ by $\chi_i(a)$. It also lies within an eigenspace of $n^{-1} a n$ for all $a \in A$, as $n^{-1} a n$ lies in $D$. It follows that $nV_{i}$ lies within an eigenspace of every $a \in A$, so it lies entirely within some $V_j$. To see this last claim, assume instead that $nV_i$ projects onto two spaces $V_{j_1}$ and $V_{j_2}$ nontrivially, and find some $b \in A$ for which $\chi_{j_1}(b) \ne \chi_{j_2}(b)$. This causes a contradiction, since $V_{j_1}$ and $V_{j_2}$ have distinct eigenvalues under $b$, while $n V_i$ lies within an eigenspace of $b$. 
\\ \\
We conclude that $n V_{i}$ is contained in some $V_j$, and conversely that $n^{-1} V_j$ projects onto $V_i$ nontrivially. By applying the above argument to $n^{-1} V_j$, we see that this is only possible if $n^{-1} V_j \subseteq V_i$, so $\dim V_i = \dim V_j$ and $n V_i = V_j$. Hence $n$ takes $V_{i}$ to $V_{j}$ for some $j$, and therefore permutes the components of the decomposition $\fv$.  This permutation defines our homomorphism $q: N \lra P$.  By definition of $q$, its kernel is precisely the stabilizer $S$.  The subgroup $S$ is in the centralizer of $D$ as scalar multiplication on $V_{i}$ commutes with every map of $V_{i}$.\\
\\
We claim that $q$ is split.  To prove this, we construct bases for each $V_{i}$ denoted by $\{v_{i,k}\}_{k=1}^{k=d_{i}}$.  For each permutation $\sigma \in P$, we know that $\sigma(i) = j$ implies $d_{i} = d_{j}$. We define a corresponding $p_{\sigma} \in N$ by defining its action on the bases.  For each $1 \leq i \leq m$, and $1 \leq k \leq d_{i}$, let $p_{\sigma}(v_{i,k}) = v_{\sigma(i),k}$.  In words, $p_{\sigma}$ takes the $k$-th basis vector of $V_{i}$ to the $k$-th basis vector of $V_{\sigma(i)}$ where $\dim V_{i} = \dim V_{\sigma(i)}$.  Extending by linearity gives our desired $p_{\sigma} \in N$.  It is routine to verify the map sending $\sigma \in P$ to $p_{\sigma} \in N$ is a group homomorphism which splits $q: N \lra P$.  
\end{proof}
The fact that the short exact sequence in Equation \ref{eq:cdn_split} is split means that each element of $N_{\GL(V)} \Diag V(\fx)$ decomposes uniquely into an element of $\Stab\left(\fv\right)$ and $\Perm( \fd )$ where we identify the latter as the subgroup of $N$ generated by the $p_{\sigma}$'s which were constructed in the proof above.  We note that both the result, and proof, of Lemma \ref{lem:complex_diag_norm} are the same if $\Diag(\fx)$ is replaced by $\Diag(\fv)$, the group of automorphisms of $V$ which act by scalar multiplication on each $V_{i}$.   We prove two more lemmas before our main theorem.

\begin{lemma} \label{lem:perm_rep_1}
Let $\rho : G \rightarrow \GL(V)$ be a faithful representation of a finite group $G = AH$ where the splitting homomorphism $\phi: H \rightarrow \Aut(A)$ is injective. There is a faithful representation $\rho'$ of $G$ that acts by generalized permutations and has the same dimension as $\rho$.
\end{lemma}
\begin{proof}
Consider the restriction of $\rho$ to $A$ and let $V = V_1 \oplus \ldots \oplus V_m$, $\fx = \{ \chi_1, \dots, \chi_m \}$ be the quantities given by Lemma \ref{lem:normalform}. For any $h \in H$, $\rho(h)$ lies in the normalizer of $\Diag V(\fx) = \rho(A)$. By Lemma \ref{lem:complex_diag_norm}, there is then a unique decomposition $\rho(h) = s(h) p(h)$, where $s(h) \in S$ and $p(h) \in P$. 
\\ \\
The map $h \mapsto p(h)$ is a representation of $H$, and the pair $(\rho|_A, p)$ further satisfies the compatibility condition of Lemma \ref{lem:semi_rep}. For any $a \in A$ and $h \in H$, 
\begin{equation*}
(\rho \circ \phi_h)(a) = \rho(h) \rho(a) \rho(h)^{-1} = s(h) p(h) \rho(a) p(h)^{-1} s(h)^{-1} = p(h) \rho(a) p(h)^{-1}.
\end{equation*}
The last equality follows because $p(h)$ and $s(h)$ lie in the normalizer and centralizer of $\rho(A)$ respectively.
\\ \\
Finally, the images of $\rho|_A$ and $p$ intersect trivially, as these lie in $S$ and $P$ respectively. As $\rho|_A$ is faithful, Lemma \ref{lem:avoids_central} ensures that $p$ is faithful, and it furnishes a faithful representation $\rho'$ of $G$ determined by $(\rho|_A, p)$. It is straightforward to check that $\rho'$ acts by generalized permutations and that $\dim \rho = \dim \rho'$, as we want.
\end{proof}
\begin{lemma} \label{lem:perm_rep_2}
Let the representations $\rho$, $\rho'$, the space $V = V_{1} \oplus \ldots \oplus V_m$, and the collection $\fx = \{ \chi_1, \dots, \chi_m \}$ be as in Lemma \ref{lem:perm_rep_1}. Then there exists a faithful representation $\rho''$ of $G$ which acts by generalized permutations and is of dimension $m$.
\end{lemma}
\begin{proof}
Following Section 8.2 of Serre \cite{serre1977linear}, we construct $\rho''$ as a faithful subrepresentation of $\rho'$. As $\rho'|_H$ acts on $V$ by permutations, for each $1 \le i \le m$, there exists a basis $\{ v_{i, k} \}_{k = 1}^{k = d_i}$ which $H$ permutes by taking $v_{i, k}$ to $v_{\sigma(i), k}$ as in the proof of Lemma \ref{lem:complex_diag_norm}. The subspace $V''$ generated by the vectors $\{v_{1, 1}, v_{2, 1}, \dots, v_{m, 1} \}$ is invariant under $P$. This subspace is also clearly invariant under the action of $D$, and thus determines a subrepresentation $\rho'': G \rightarrow \GL(V'')$ of dimension $m$ which acts by generalized permutations. This subrepresentation is faithful, as $\rho''(g)$ fixes each $v_{i, 1}$ only if $\rho'(g)$ fixes each $V_i$.
\end{proof}
In the proofs of Lemmas \ref{lem:complex_diag_norm}, \ref{lem:perm_rep_1}, and \ref{lem:perm_rep_2} care is taken to not assume that each $\dim V_{i} = 1$.  This is typically not true, so one must work with the isotypic components of the representation of $A$ to construct the permutation representation.  An earlier preprint of this work missed this subtlety.   Lemma \ref{lem:perm_rep_2} shows if we are interested in minimizing the dimension, we may always work under this assumption.  

\begin{theorem}\label{thm:cmdim}
The minimal dimension of a faithful representation of $G = AH$ with faithful splitting homomorphism is equal to the smallest size of an $H$-invariant generating set of $\hat{A}$.  Moreover, $G$ admits a minimal faithful representation into the generalized permutation group.  
\end{theorem}
\begin{proof}
The latter claim follows immediately from Lemma \ref{lem:perm_rep_2} because any faithful representation may be modified to take its values in the generalized permutation group.  In particular, this applies to any given minimal faithful representation.  We therefore need to only prove the dimension claim.  
\\
\\
Let $\fx = \{\chi_{1},\hdots, \chi_{m}\}$ be a smallest size $H$-invariant generating set of $\hat{A}$.  We will use $\fx$ to construct a faithful representation of $G$.  Begin by forming the vector space generated by the basis $\{e_{\chi_{1}},\hdots, e_{\chi_{m}}\}$ so that $V = V_{1}\oplus \hdots \oplus V_{m}$ where $V_{i} := \C e_{\chi_{i}} $ as was done in the proof of Theorem \ref{thm:min_abl}.  We define a representation of $G$ using Lemma \ref{lem:semi_rep}.  For each $a \in A$, let $\rho(a)$ act on $V_{i}$ by $\chi_{i}$.  For each $h \in H$ define $\sigma(h): V \lra V$ on the basis for $V$ via $\sigma(h)(e_{\chi_{i}}) := e_{h\chi_{i}}$ where $h\chi_{i}$ is given by the permutation action of $h$ on $\fx$.  By construction, $\rho: A \lra \GL(V)$ takes its values in $D$ and $\sigma$ takes its values in $P$.  We now only need to verify the compatibility condition of Equation \ref{eq:comp-rel}.  \\
\\
Let $a \in A$, $h \in H$, and pick an $e_{\chi_{i}}$.  Then 
\begin{align} \label{eq:chiinv}
&\left(\sigma(h)\rho(a)  \sigma(h)^{-1}\right)(e_{\chi_{i}}) = \left(\sigma(h) \rho(a)\right)e_{h^{-1}\chi_{i}} = \sigma(h)\left(\left(h^{-1}\chi_{i}\right)(a)e_{h^{-1}\chi_{i}}\right)\\
&= (h^{-1}\chi_{i})(a)e_{\chi_{i}} = \chi_{i}(hah^{-1})e_{\chi_{i}} =  \rho(hah^{-1})e_{\chi_{i}} = \left(\rho\circ\phi_{h}\right)(a) (e_{\chi_{i}}) \notag
\end{align}
As the equality holds on a basis of $V$, it holds on all of $V$ thus verifying the desired compatibility condition, hence $\rho$ extends to a representation of all of $G = AH$.  We claim this extension is faithful.  Because $\fx$ is a generating set, $\rho$ on $A$ is faithful by Lemma \ref{lem:dual_group_gen}.  Because the image of $\sigma$ intersects the image of $\rho$ trivially, faithfulness of the whole representation is equivalent to $\sigma$ on $H$ being faithful by Lemma \ref{lem:avoids_central}.  The fact that $\sigma$ on $H$ is faithful is due to Lemma \ref{lem:dual_action_faithful}, as $\fx$ is a generating set of $\hat{A}$, so $\rho: G \lra \GL(V)$ is faithful.  \\
\\
Thus given a minimal $H$-invariant generating set of $\hat{A}$ of size $m$, we may construct a faithful representation of $G$ of dimension $m$.  To conclude the proof, we need to show the dimension of any minimal faithful representation of $G$ is bounded below by a minimal $H$-invariant generating set of $\hat{A}$.  Given a minimal faithful representation of $G$, Lemma \ref{lem:perm_rep_2} guarantees that the dimension is bounded below by the number of distinct characters $\fx$ associated to the representation of $A$.  By Lemma \ref{lem:dual_group_gen}, this set generates $\hat{A}$, and by construction, the collection is $H$-invariant.  Thus, the dimension of any minimal faithful representation of $G$ is bounded below by an $H$-invariant generating set for $\hat{A}$, and in particular must be at least as large as the smallest one.  
\end{proof}
We observe that Theorem \ref{thm:cmdim} justifies our previous claim that in some sense groups $G = AH$ with injective splitting homomorphism are generalized permutation groups as there exist minimal faithful representations of them into this group.  It follows from Theorem \ref{thm:cmdim} that $\mdimc G $ is thus equal to the dimension of the smallest torus which we may embed $A$ into in such a manner that the action of $H$ on $A$ is represented by symmetries of the torus.  This is made precise in the corollary below.

\begin{corollary}\label{cor:toremb}
The minimal dimension of a faithful representation of $G = AH$ with an injective splitting homomorphism is equal to the smallest dimension $m$ for which there exists both an embedding of $\psi: A \lra \Z_{d}^{m}$ for some $d \geq 2$ and an injective homomorphism $\sigma: H \lra S_{m}$ which is equivariant with respect to the standard permutation action of $S_{m}$ on $\Z_{d}^{m}$.  Symbolically, the diagram below commutes for each $h \in H$.  
\begin{equation}\label{eq:toruseq}
\begin{tikzcd}
A \arrow[r, "\psi"] \arrow{d}[swap]{\phi_{h}}
& \Z_{d}^{m} \arrow[d, "\sigma_{h}"] \\
A \arrow{r}[swap]{\psi}
& \Z_{d}^{m}
\end{tikzcd}
\end{equation}
\end{corollary}

\begin{proof}
Choose a minimal faithful representation of $G$ of dimension $m$ that acts on $V = V_{1}\oplus \hdots \oplus V_{m}$ by generalized permutations with distinct characters $\fx = \{\chi_{1},\hdots, \chi_{m}\}$.  As each element of $\fx$ has at most order $d := |A|$, for each $a \in A$, we may identify $\chi_{i}(a)$ with an element of $\Z_{d}$, and thus identify $a$ with the tuple $(\chi_{1}(a),\hdots , \chi_{m}(a)) \in \Z_{d}^{m}$.  This determines our embedding $\psi: A \lra \Z_{d}^{m}$.  In addition, each element $h \in H$ acts on $V = V_{1}\oplus \hdots \oplus V_{m}$ by permutations which respect the decomposition.  This determines the action of $H$ on $\Z_{d}^{m}$ by the standard permutation action.  Hence every faithful representation of $G$ determines a pair of embeddings $\psi: A \lra \Z_{d}^{m}$ and $\sigma: H \lra S_{m}$ which satisfy Equation \ref{eq:toruseq}.  Thus $\mdimc G $ is bounded below by the minimal dimension of $\Z_{d}^{m}$ satisfying the condition in the hypotheses of the corollary.       \\
\\
Next let $m$ be the smallest such dimension satisfying the condition in the hypotheses of the corollary and $\psi: A \lra \Z_{d}^{m}$ and $\sigma: H \lra S_{m}$ be the associated pair of embeddings.  It is routine to verify the collection of characters of $A$ defined by $\chi_{i} := p_{i}\circ \psi: A \lra S^{1}$ where $p_{i}: \Z_{d}^{m} \lra \Z_{d} \leq S^{1}$ is the projection onto the $i$-th coordinate of $\Z_{d}^{m}$ is an $H$-invariant generating set of $\hat{A}$.  The fact that $\psi: A \lra \Z_{d}^{m}$ is an embedding shows that the collection $\fx$ generates $\hat{A}$, and $H$-invariance follows from equivariance of the embedding with respect to the actions of $H$ on $A$ and $\Z_{d}^{m}$.  By Theorem \ref{thm:cmdim}, this means $m \geq \mdimc G $, thus completing the proof.
\end{proof}

\section{Real Representations of Split Extensions by Abelian Groups}\label{sec:rrseag}
In this section we prove analogous results from Section \ref{sec:rseag} in the context of real representations of groups of the form $G = AH$ with faithful splitting homomorphism $\phi: H \lra \Aut(A)$.  Most of the results from Section \ref{sec:rseag} derive in some manner from Lemma \ref{lem:complex_diag_norm}.  Our first goal is to prove a real analogue of this lemma.  \\
\\
We begin with notation that will persist throughout this section.  Let $Z = U\oplus W = \left(U_{1}\oplus \hdots \oplus U_{n}\right) \oplus \left(W_{1} \oplus \hdots \oplus W_{m}\right)$ where $\dim U_{i} = c_{i}$, $\dim W_{i} = 2d_{i}$, and let $\fc$, $\fd$ be the ordered list of dimensions of $U_{i}$ and $W_{i}$ respectively.  Let $J$ and $I$ be compatible real and almost complex structures on $Z$ and $W$ which both respect the decompositions $\fz$ and $\fw$ respectively.  Let $\mathcal{R} \cup \mathcal{X} : =\{\eta_{1},\hdots, \eta_{n}\} \cup \{\chi_{1},\hdots, \chi_{m}\}$ be a collection of distinct real characters and distinct pairwise non-conjugate complex characters of $A$.  We will denote $\fx^{\pm} :=  \{\chi_{1}, \ol{\chi_{1}}, \hdots, \chi_{m},\ol{\chi_{m}}\}$ to be the collection of distinct complex characters $\fx$ with their conjugates included.  \\
\\
Let $D:=\Diag Z(\mathcal{R},\mathcal{X})$ be the subgroup of real automorphisms $\GL(Z,J)$ which acts on each $U_{i}$ by $\eta_{i}$ and $(W_{i},J,I)$ by $\chi_{i}$.  Let $S:=\Stab(\fz,J,I)$ be the subgroup of $\GL(Z,J)$ which preserves the decomposition $\fz$, the real structure $J$, and the almost-complex structure $I$ on $W$.  Let $\Perm(\fc,\fd^{\pm})$ be the permutations of the set $\{1_{U},\hdots, n_{U}\} \cup \{1_{W}^{+},1_{W}^{-},\hdots, m_{W}^{+}, m_{W}^{-}\}$ which permute the indices $i_{U},j_{U}$ of $U$ with the same dimensions $c_{i} = c_{j}$ and the pairs of indices $\{i_{W}^{+},i_{W}^{-}\}, \{j_{W}^{+},j_{W}^{-}\}$ of $W$ with the same dimensions $2d_{i} = 2d_{j}$. Finally, let $\Perm(\fc,\fd)$ be the subgroup of $\Perm(\fc,\fd^{\pm})$ which take positive $W$-indices to positive $W$-indices and negative $W$-indices to negative $W$-indices.
%

\begin{lemma}\label{lem:real_diag_norm}
Using the notation in the above paragraphs, there exists a surjective homomorphism $q$ from the normalizer of $\Diag Z(\fr,\fx)$ inside of $\GL(Z,J)$ onto $\Perm(\fc,\fd^{\pm}) \leq S_{n}\times \osym_{m}$ whose kernel is $\Stab(\fz, J, I)$ which centralizes $\Diag Z(\fr,\fx)$.  Moreover, this short exact sequence is split.  Symbolically,
\begin{equation}\label{eq:gen_cdn_split}
\Stab(\fz,J,I) \xrightarrow{i} N_{\GL(Z,J)} \Diag Z(\fr,\fx) \xrightarrow{q} \Perm(\fc,\fd^{\pm})
\end{equation}
is split.
\end{lemma}

\begin{proof}
We let $S, D, N$ and $P$ denote the obvious groups appearing in Equation \ref{eq:gen_cdn_split}.  Observe the action of $H$ on the dual of $A$ sends real characters to real characters and complex characters to complex characters.  Consequently the $H$-action must send $U$ to $U$ and $W$ to $W$.  We claim this action in addition preserves the decompositions $\fu$ of $U$ and $\fw$ of $W$, which will induce our homomorphism $q$. 
\\ \\
We prove this in two steps, starting with $\fu$. For any $n \in N$ and $U_i \in \fu$, we may recycle the proof of Lemma \ref{lem:complex_diag_norm} to show that $n(U_i) = U_j$ for some $U_j \in \fu$. Thus to each $n \in N$ we assign the permutation $\sigma(n)(i_U) = j_U$ if $n(U_i) = U_j$.
\\ \\
We now turn to $\fw$. By definition, for each $1 \le k \le m$, $a \in A$ acts by the automorphism $\re \chi_k(a) \, \id_{W_k} + \im \chi_k(a) \, I$ on each $W_k$. From the remarks after Definition \ref{def:actsbychi}, we recall that $W_k$ decomposes into spaces $W_k^+$ and $W_k^-$ with eigenvalues of $i$ and $-i$ under $I$. The spaces $W_k^\pm$ lie in eigenspaces of each $a \in A$, as $a$ acts on $W_k^+$ by $\chi_k(a)$ and on $W_k^-$ by $\overline{\chi_k}(a)$. As the $\chi_k$ are complex, distinct, and pairwise non-conjugate, any two of these spaces must also have distinct eigenvalues under some $a \in A$. It then follows from the proof of Lemma \ref{lem:complex_diag_norm} that any $n \in N$ permutes the spaces $W_k^{\pm}$.
\\ \\
To see that this permutation is oriented, take any space $W_k$ with subspaces $W_k^+$, $W_k^-$. From Definition \ref{def:actsbychi}, we further recall that $J$ interchanges these subspaces. Choose any $n \in N$, and let $n(W_k^+) = W_{l}^{\pm}$ for some $1 \le l \le m$. As $n$ lies in $\GL(Z, J)$, we find that
\begin{equation*}
n (W_k^-) = n J(W_{k}^+) = Jn(W_k^+) = J(W_{l}^{\pm}) = W_{l}^{\mp}.
\end{equation*}
We then obtain an oriented permutation $\tau(n)$ of the set $\{ k_{W}^{\pm} \}$ by defining $\tau(n)(k_W^+) = l_W^\pm$ and $\tau(n)(k_W^-) = l_W^{\mp}$ if $n(W_k^+) = W_l^\pm$. We combine the permutations $\sigma$ and $\tau$ as defined above to obtain our homomorphism $q: N \lra P$ defined on all of $\{1_{U},\hdots, n_{U}\} \cup \{1_{W}^{+},1_{W}^{-},\hdots, m_{W}^{+},m_{W}^{-}\}$.  It is a routine, albeit tedious, calculation to verify $q$ is a homomorphism.  
\\
\\
Elements $n \in N$ in the kernel of $q$ by definition satisfy the property that they stabilize the entire decomposition of $\fz$ as they stabilize each $U_{i}$ and $W_{i}$.  Thus we have the following equalities below for $1 \leq i \leq m$ and any $a \in A$ restricted to $W_{i}$.  
\begin{align}\label{eq:n_preserves_i}
nan^{-1}|_{W_{i}} &= n\left(\re \chi_{i}(a)\id + \im \chi_{i}(a) I\right) n^{-1} = \re \chi_{i}(a) \id + \im \chi_{i}(a) (nIn^{-1}) \notag \\
a'|_{W_{i}} &=\re \chi_{i}(a')\id + \im \chi_{i}(a') I
\end{align}
If we choose an $a \in A$ where $\im \chi_{i}(a) \neq 0$, we may restrict the maps above in Equation \ref{eq:n_preserves_i} to $W_i^+$, and equate real and imaginary parts of the rightmost quantities.  In so doing, we see that $\im \chi_{i}(a) nIn^{-1} = \im \chi_{i}(a') I$ as maps from $W_{i}$ to itself.  Thus, we may find $\lambda_{i}$'s so that $nIn^{-1} = \lambda_{i} I$ for all $1 \leq i \leq m$.  Because the eigenvalues of $I$ are $\pm i$, this means $\lambda_{i} = \pm 1$.  Because $q(n) = 1$, $n(W_i^+) = W_i^+$ so $\lambda_{i} = 1$ for all $1 \leq i \leq m$.  In particular, this means that if $q(n) = 1$ then $n$ is contained in $S$.  That $S$ is both in the kernel of $q$, and, in the centralizer of $D$ follows readily from Equation \ref{eq:n_preserves_i}.  
\\
\\
Finally, we prove this homomorphism is split.  We construct bases for each $U_{i}$ and $W_{i}$ adapted to the structures $J$ and $I$.  For $1 \leq i \leq n$, let $\{U_{i,k}\}_{k=1}^{k=c_{i}}$ be a basis for $U_{i}$.  For $1 \leq i \leq m$, let $\{w_{i,k}\}_{k=1}^{k=d_{i}}$ be a basis for $W_i^+$ and $\{\ol{w_{i,k}}\}_{k=1}^{k=d_{i}}$ denote the images of each $w_{i,k}$ in $W_i^-$ under $J$.  We call the collections of $w_{i,k}$'s and $\ol{w_{i,k}}$'s positive and negative, respectively.  For each permutation $\sigma \in P$, we define a corresponding $p_{\sigma} \in N$ by defining its actions on the bases.  To make notation less cumbersome, we drop the subscripts $U$ and $W$ on the indices $1 \leq i_{U} \leq n$ and $1 \leq i_{W}^{\pm} \leq m$.  \\
\\
For each $1 \leq i \leq n$ and $1 \leq k \leq c_{i}$, let $p_{\sigma}(u_{i,k}) := u_{\sigma(i),k}$.  In words, $p_{\sigma}$ takes the $k$-th basis vector of $U_{i}$ to the $k$-th basis vector of $U_{\sigma(i)}$.  For each $1 \leq i^{+} \leq m$ and $1 \leq k \leq d_{i}$, let $p_{\sigma}(w_{i,k}) := w_{\sigma(i),k}$ if $\sigma(i^{+})$ is positive, or, $p_{\sigma}(w_{i,k}) := \ol{w_{\sigma(i),k}}$ if $\sigma(i^{+})$ is negative.  Define $p_{\sigma}(\ol{w_{i,k}})$ similarly.  In words, $p_{\sigma}$ takes the $k$-th positive basis vector of $W_{i}$ to the $k$-th positive, or negative, basis vector of $W_{\sigma(i)}$ depending on the sign of $\sigma(i^{+})$.  A similar description is true for $p_{\sigma}$ on the $k$-th negative basis vector of $W_{i}$.  By construction, the assignment taking $\sigma \in P$ to $p_{\sigma} \in N$ is a homomorphism from $P$ to $N$ which splits $q: N \lra P$.  
\end{proof}

We remark that Lemma \ref{lem:real_diag_norm} is a structural statement about the normalizer of $\Diag Z(\mathcal{R},\mathcal{X})$ inside $\GL(Z,J)$.  The subgroup $\Stab(Z,J,I)$ lies inside the normalizer as a finite-index subgroup, thus `most' elements of the normalizer of $\Diag Z(\mathcal{R},\mathcal{X})$ preserve the almost-complex structure $I$ on $W$.  In fact, the obstruction to an element preserving the almost-complex structure $I$ is measured precisely by how it decomposes in $\Perm(\fc,\fd^{\pm})$.  We make this statement precise in the following corollary.  

\begin{corollary}\label{cor:almost_preserve}
An element in $N_{\GL(Z,J)} \Diag Z(\fr,\fx)$ preserves the almost complex structure $I$ on $W$ if and only if its image under $q$ as defined in Lemma \ref{lem:real_diag_norm} takes its values inside $\Perm(\fc, \fd) \leq \Perm(\fc, \fd^{\pm})$.  
\end{corollary}

\begin{proof}
We adopt the same notation as in the proof of Lemma \ref{lem:real_diag_norm}.  As seen in the proof of this lemma, $n \in N$ must permute the spaces $W_{k}^{\pm}$ for each $1 \leq k \leq m$.  If $n$ in addition preserves the almost-complex structure $I$, then $n$ must send $W_{k}^{+}$ to $W_{l}^{+}$ and $W_{k}^{-}$ to $W_{l}^{-}$.  Consequently, the corresponding permutation $q(n)$ as defined on $\{1_{U},\hdots, n_{U}\} \cup \{1_{W}^{+},1_{W}^{-},\hdots, m_{W}^{+},m_{W}^{-}\}$ takes its values inside $\Perm(\fc,\fd)$ as claimed.  \\
\\
To prove the converse direction, observe that if we let $\sigma := q(n)$, then $p_{\sigma} \in N$ as defined in Lemma \ref{lem:real_diag_norm} for $\sigma \in \Perm(\fc,\fd)$ preserves the almost-complex structure $I$.  By the conclusion of Lemma \ref{lem:real_diag_norm}, $n$ thus decomposes into the composition of an element in $\Stab(\fz, J, I)$ and $p_{\sigma}$, both of which preserve $I$.  
\end{proof}

Before proving the main theorem of this section, we state a real analogue of Lemma \ref{lem:perm_rep_2} that allows us to modify a given faithful real representation.  The proof follows from the same arguments in Lemmas \ref{lem:perm_rep_1} and \ref{lem:perm_rep_2} using the machinery developed in Section \ref{sec:preliminaries} and Lemma \ref{lem:real_diag_norm} and is thus omitted.  
\begin{lemma}\label{lem:real_perm_rep}
Let $\rho: G \lra \GL(Z,J)$ be a faithful representation of a finite group $G = AH$ where the splitting homomorphism $H \lra \Aut(A)$ is injective.  Let $\fr \cup \fx^{\pm} = \{\eta_{1},\hdots, \eta_{n}\} \cup \{ \chi_{1},\ol{\chi_{1}},\hdots, \chi_{m}, \ol{\chi_{m}}\}$ be the collection of distinct characters of the representation of $A$ on $Z$. Then there exists a faithful real representation $\rho$ of $G$ acting by generalized permutations and whose dimension is equal to $n + 2m$.  
\end{lemma}
 
\begin{theorem}\label{thm:rmdim}
The minimal dimension of a faithful real representation of $G = AH$ with an injective splitting homomorphism is equal to the smallest size of an $H$-invariant and conjugation-invariant generating set of $\hat{A}$.  Moreover, $G$ admits a minimal faithful representation into the generalized permutation group.   
\end{theorem}

\begin{proof}
The proof of the latter statement follows from Lemma \ref{lem:real_perm_rep} and the arguments in the beginning of Theorem \ref{thm:cmdim}.  If we are given a smallest size $H$-invariant and conjugation-invariant generating set of $\hat{A}$, we partition it into real and complex characters say, $\mathcal{R}:=\{\eta_{1},\hdots, \eta_{n}\}$ and $\mathcal{X}^{\pm} := \{\chi_{1},\ol{\chi_{1}},\hdots, \chi_{m},\ol{\chi_{m}}\}$.  We construct a vector space $Z$ generated by the basis $\{e_{\eta_{1}},\hdots, e_{\eta_{n}}\}  \cup \{e_{\chi_{1}}, e_{\ol{\chi_{1}}}, \hdots, e_{\chi_{m}}, e_{\ol{\chi_{m}}}\}$ so that $Z = U\oplus W :=(U_{1} \oplus \hdots \oplus U_{n}) \oplus \left(W_{1} \oplus \hdots \oplus W_{m} \right)$ where $U$ is the sum of the $U_{i}:= \C e_{\eta_{i}}$, and $W$ is the sum of the $W_{i} := \C e_{\chi_{i}} \oplus \C e_{\ol{\chi_{i}}}$.  Define a real structure $J$ on $Z$ which takes $e_{\eta_{i}}$ to $e_{\eta_{i}}$ for $1 \leq i \leq n$ and swaps $e_{\chi_{i}}$ and $e_{\ol{\chi_{i}}}$ for $1 \leq i \leq m$, then extend \emph{conjugate-linearly} to all of $Z$.  Finally, define an almost-complex structure $I$ on $W$ by multiplication by $i$ on the space spanned by $\{e_{\chi_{i}}\}_{i=1}^{i=m}$ and by multiplication by $-i$ on the space spanned by $\{e_{\ol{\chi_{i}}}\}_{i=1}^{i=m}$.   \\
\\
Define a representation $\rho:A \lra \Diag Z(\mathcal{R},\mathcal{X})$ by acting on each $U_{i}$ by $\eta_{i}$ and each $(W_{i},J,I)$ by $\chi_{i}$.  Define $\sigma: H \lra N_{\GL(Z,J)} \Diag Z(\mathcal{R},\mathcal{X})$ via $\sigma(h)(e_{\nu_{i}}) := e_{h\nu_{i}}$ for each $\nu \in \mathcal{R}\cup \mathcal{X}^{\pm}$.  The same calculation as in Equation \ref{eq:chiinv} shows the pair $(\rho,\sigma)$ defines a compatible pair in the sense of Lemma \ref{lem:semi_rep}, and thus defines a representation of $G = AH$.   Because the image of $\sigma$ intersects the image of $\rho$ trivially, Lemma \ref{lem:avoids_central} guarantees faithfulness of this representation if we know both $\rho$ and $\sigma$ are faithful.  Both are faithful by Lemmas \ref{lem:dual_group_gen} and \ref{lem:dual_action_faithful} because $\fr\cup \fx^{\pm}$ generates $\hat{A}$.   
 \\
\\
To conclude the proof, we need to show the dimension of any minimal faithful real representation of $G = AH$ is bounded below by the size of a smallest $H$-invariant and conjugate invariant generating set of $\hat{A}$.  This follows from the same arguments at the end of Theorem \ref{thm:cmdim} and Lemma \ref{lem:real_perm_rep}.
\end{proof}

Theorem \ref{thm:rmdim} both addresses the minimal dimension of $G = AH$ and illuminates the structure of a such a minimal faithful real representation by furnishing a generalized permutation representation.  However, Lemma \ref{lem:real_perm_rep} shows that up to a finite index, the representation is an almost-complex representation.  Admitting a minimal faithful almost-complex representation imposes a strong restriction on the representation structurally, but it can be entirely detected by group invariants as seen in the corollary below.  
\begin{corollary}\label{cor:when_ac}
The group $G = AH$ with a faithful splitting homomorphism admits an almost-complex minimal faithful representation if and only if $\hat{A}$ admits a smallest size $H$-invariant and conjugation invariant generating set consisting entirely of complex characters, $\fx^{\pm} =  \{ \chi_{1},\ol{\chi_{1}},\hdots, \chi_{m}, \ol{\chi_{m}}\}$, where the $H$-action on $\fx^{\pm}$ preserves the decomposition of conjugate and non-conjugate complex characters.  
\end{corollary}

\begin{proof}
If $\hat{A}$ admits a smallest size $H$-invariant conjugation invariant generating set where the $H$-action preserves the decomposition, Theorem \ref{thm:rmdim} guarantees the existence of an almost-complex minimal faithful representation as the almost-complex structure $I$ as constructed in the proof is invariant under $H$. Indeed, the $H$-action preserves the decomposition $\{e_{\chi_{i}}\}_{i=1}^{m}$ and $\{e_{\ol{\chi_{i}}}\}_{i=1}^{m}$ which are precisely the $\pm i$-eigenspaces of $I$.  \\
\\
Conversely, let $\rho: G \lra \GL(Z,J,I)$ be an almost-complex minimal faithful representation.  Recalling Lemma \ref{lem:real_perm_rep}, we first show that the characters $\fr \cup \fx^{\pm}$ of $\hat{A}$ are all complex.  To this end, we restrict the representation to $A$ and obtain the decomposition $Z = U \oplus W$ where $A$ acts on $U$ by real characters and on $W$ by complex characters.  Because $A$ preserves both structures $J$ and $I$, $U$ is invariant under both $J$ and $I$.  Thus, we obtain an almost-complex representation $\rho: A \lra \GL(U,J,I)$.  By Proposition \ref{lem:almostcomplexeq}, this means $(U,J) \cong (V\oplus \ol{V}, J')$ for some complex representation $V$ of $A$.  However, because $A$ acts on $U$ by entirely real characters, this means $V \cong \ol{V}$. This contradicts minimality by Lemma \ref{lem:real_perm_rep}, as the characters of $A$ that arise from a minimal real representation of $G$ must be distinct.  Hence $U = 0$.  Thus, the corresponding collection of characters is entirely complex.  That this is the smallest size $H$-invariant and conjugation-invariant collection follows from Theorem \ref{thm:rmdim}, and the fact it preserves a decomposition follows from Corollary \ref{cor:almost_preserve}.
\end{proof}

Note as a consequence of Corollary \ref{cor:when_ac} and Theorem \ref{thm:min_abl} it follows that a finite abelian group admits an almost-complex minimal faithful real representation if and only if the invariant factors of $A$ are all greater than $2$.  We also remark that in the case of Corollary \ref{cor:when_ac}, we have that $\mdimr G = 2\mdimc G$.  This follows more generally for groups that admit almost-complex minimal faithful representations.  In this case, the representation decomposes as $V\oplus \ol{V}$ for some complex representation by Lemma \ref{lem:almostcomplexeq}.  Because the kernels of $V$ and $\ol{V}$ are equal, and the kernel of the sum is trivial by hypothesis, $V$ is a faithful complex representation of $G$.  To see this complex dimension is minimal, observe if there were a smaller faithful complex representation $V'$, then $V' \oplus \ol{V'}$ would be a smaller real representation, again by Lemma \ref{lem:almostcomplexeq}.
\\
\\
We conclude this section by stating the real analogue to Corollary \ref{cor:toremb}.  It says that the minimal faithful real dimension of $G = AH$ is the smallest dimension in which we can represent the $H$ action on $A$ by a combination of permutations and oriented permutations.  The proof is omitted as it is nearly identical to the proof of Corollary \ref{cor:toremb}.  
\begin{corollary}\label{cor:real_toremb}
The minimal dimension of a real faithful representation of $G = AH$ with an injective splitting homomorphism is equal to the smallest dimension $n+2m$ for which there exists both an embedding of $\psi: A \lra \Z_{2}^{n}\times \Z_{d}^{2m}$ for some $d > 2$ and an injective homomorphism $\sigma: H \lra S_{n}\times S_{m}^{\pm}$ which is equivariant with respect to the action of $S_{n}\times S_{m}^{\pm}$ on $\Z_{2}^{n}\times \Z_{d}^{2m}$ where $S_{n}$ acts on the $\Z_{2}^{n}$ by the standard permutation action and $S_{m}^{\pm}$ acts on $\Z_{d}^{2m}$ by oriented permutations.  Symbolically, the diagram below commutes for each $h \in H$.  
\begin{equation}\label{eq:real_toruseq}
\begin{tikzcd}
A \arrow[r, "\psi"] \arrow{d}[swap]{\phi_{h}}
& \Z_{2}^{n}\times \Z_{d}^{2m} \arrow[d, "\sigma_{h}"] \\
A \arrow{r}[swap]{\psi}
& \Z_{2}^{n}\times \Z_{d}^{2m}
\end{tikzcd}
\end{equation}
\end{corollary}


\section{Applications}\label{sec:applications}
We conclude the paper with some applications.  In Proposition \ref{prop:rubiks}, we compute the minimal real and complex dimensions of the Rubik's cube group described in Example \ref{ex:rubiks}. We then introduce a subfamily of groups for which Moret\'o's bound is sharp and exhibit others for which $\mdimc G $ is logarithmic in $|G|$. 
 In addition, we show that one cannot always obtain a minimal real representation by adding conjugate pieces to some minimal complex representation.  We conclude by studying some of the structure of these minimal faithful real representations and their corresponding orbifold singularities.  

\begin{proposition} \label{prop:rubiks}
Let $R_{3}$ denote the $3\times 3$ Rubik's cube group.  The minimal complex and real faithful dimensions of $R_{3}$ are 20 and 28 respectively.  
\end{proposition}
\begin{proof}
As noted in Example \ref{ex:rubiks}, $R_{3}$ is isomorphic to $AH$ where $A = \Z_{2,0}^{12}\oplus \Z_{3,0}^{8}$ and $H$ is the subgroup of permutations $(\alpha,\beta) \in S_{12}\times S_{8}$ with $\sgn(\alpha) = \sgn(\beta)$.  By Corollary \ref{cor:toremb}, $H$ embeds into $S_{m}$ where $m = \mdimc R_{3} $, thus $A_{8}\times A_{12}$ embeds into $S_{m}$.  This necessitates that the size of a minimal $H$-invariant generating set of $A$ is bounded below by $8 + 12 = 20$.  \cite{wright}  \\
\\
Thus to prove $\mdimc R_{3}  = 20$, it suffices to exhibit an $H$-invariant generating set of the dual of $A:= \Z_{2,0}^{12}\oplus \Z_{3,0}^{8}$ with $20$ elements.  Consider the collection $\mathcal{R} \cup \mathcal{X} := \{\eta_{1},\hdots, \eta_{12}\} \cup \{\chi_{1},\hdots , \chi_{8}\}$ where $\eta_{i}: A \lra \Z_{2} \leq S^{1}$ and $\chi_{i} : A \lra \Z_{3} \leq S^{1}$ are the projections onto the $i$-th coordinates of $\Z_{2,0}^{12}$ and $\Z_{3,0}^{8}$ respectively.  It is clear this collection is $H$-invariant.  Because their kernels intersect trivially, Lemma \ref{lem:dual_group_gen} implies the collection is also generating, thus $\mdimc R_{3} = 20$.  
\\
\\
We now address the real case.  By Corollary \ref{cor:real_toremb}, the size of a minimal $H$-invariant and conjugation-invariant generating set of $A$ must admit an embedding of $H$ into $S_{n}\times \osym_{m}$.  Immediately we see that $n + m \geq 20$.  We claim that $m \geq 8$.  Consider the projection of $H \lra \osym_{m}$.  This map cannot be trivial as there must exist complex characters to represent $\Z_{3,0}^{8}$.  By considering possible kernels of the $H \lra \osym_{m}$, the image of this map inside of $\osym_{m}$ contains either $A_{8}, A_{12},$ or, $A_{8}\times A_{12}$ as $H \cong (A_{12}\times A_{8})\rtimes \Z_{2}$.  This necessitates that $m \geq 8$, thus the size of a minimal $H$-invariant and conjugation-invariant generating set of $A$ is bounded below by $28$.\\
\\
To provide such a generating set of $A$, we simply take the previous generating set $\mathcal{R}\cup \mathcal{X}$, and close it under conjugation to yield $\{\eta_{1},\hdots, \eta_{12}\} \cup \{\chi_{1},\ol{\chi_{1}},\hdots, \chi_{8}, \ol{\chi_{8}}\}$.  Consequently, $\mdimr R_{3} = 28$ as claimed.  
\end{proof}
The minimal complex dimension of $R_3$ was independently determined by D. Holt \cite{HoltRubiks}.
\\ \\
We now briefly consider whether Moret\'o's bound may be improved for our class of groups. The following example shows that for every prime $p$, there is a minimal representation of dimension $p - 1$ that saturates Moret\'o's inequality. However, many minimal representations in the example are more efficient.

\begin{example} \label{ex:affine_groups}
For any prime $p$, the affine group $\Aff(n, \Z_p)$ acts on the vector space $\Z_p^n$ by translations and linear transformations. The subgroup $A \cong \Z_p^n$ of translations is abelian and normal in $\Aff(n, \Z_p)$, and the subgroup $H \cong \GL(n, \Z_p)$ of linear transformations is complementary to $A$. The conjugation action of $H$ on $A$ is faithful, and it is transitive on the nontrivial elements of $A$. It follows that the smallest $H$-invariant generating set for the dual group $\hat{A}$ is of size $|A| - 1 = p^n - 1$. By Theorem \ref{thm:cmdim} and Theorem \ref{thm:rmdim}, this is the minimal dimension of $\Aff(n, \Z_p)$ over both $\C$ and $\R$. Meanwhile, the order of $\Aff(n, \Z_p)$ is $|A| |H| = p^n[(p^n - 1)(p^n - p) \dots (p^n - p^{n - 1})]$.
\\ \\
We then find that $\mdimc G  \approx \sqrt[n + 1]{|G|}$, which surpasses Moret\'o's bound $\mdimc G \le \sqrt{|G|}$ for large $n$. For $n = 1$, however, we find
\[ \mdimc \Aff(1, \Z_p) = p - 1 = \text{floor} \, \sqrt{p(p - 1)} = \text{floor} \, \sqrt{|\Aff(1, \Z_p)|}. \] 
This indicates that Moret\'o's inequality is sharp among our class of groups, although many members of the class have minimal dimension well below Moret\'o's bound.
\end{example}

To emphasize that Moret\'o's bound is not always optimal, we also describe a family of groups whose minimal dimension is logarithmic in the order of the group.


\begin{example} \label{ex:generalized_rubiks}
For any positive integers $n$ and $m > 1$, let $G(n, m)$ be the semi-direct product $\Z_{m}^n \rtimes S_n$, where $S_n$ acts by permuting coordinates as before.  Similarly, following the notation in Example \ref{ex:rubiks}, let $G_0(n, m) = \Z_{m, 0}^n \rtimes S_n$.  We are interested in calculating both the minimal complex and real faithful dimensions of these groups.  As both $G(1,m)$ and $G_{0}(1,m)$ are abelian, and their minimal dimensions are calculated in Theorem \ref{thm:min_abl}, we may assume that $n > 1$.  
\\
\\
By generalizing Proposition \ref{prop:rubiks}, we readily find that $\mdimc G_0(n, m) = n$. Moreover, we may construct an $n$-dimensional representation of $G(n, m)$ by allowing $S_n$ to permute the $n$ irreducible components of a minimal representation of $\Z_m^n$. Since $G_0(n, m) < G(n, m)$, this implies that $\mdimc G(n, m) = n$ as well.  Given that $|G(n, m)| = m^n n!$ and $|G_0(n, m)| = m^{n - 1} n!$, the minimal dimensions of these families are again well within Moret\'o's bound.\\
\\
We now turn to the real case. By extending the minimal real representations of $\Z_m^n$, we find that $\mdimr G(n, m) = n$ if $m = 2$ and $2n$ otherwise. We further deduce that $\mdimr G_0(n, 2) = n$, as the minimal representations of this family over $\C$ are real-valued.
\\ \\
To determine $\mdimr G_0(n, m)$ for $m > 2$, we first note that if $n \ne 2, 4$, $S_n$ embeds in $\osym_c$ only if $c \ge n$. Indeed, an embedding $S_n \longrightarrow \osym_c$ produces a map $S_n \lra \osym_c \lra S_c$ whose kernel is a 2-group. Every such map is injective except the quotients $S_2 \lra S_1$ and $S_4 \lra S_3$, with kernels $\Z_2$ and the Klein four group. In all other cases, $S_n$ embeds in $S_c$ and $c \ge n$ as we want.
\\ \\
We then claim that $\mdimr G_0(n, m) = 2n$ for $m > 2$ and $n \ne 2, 4$. To see this, recall the embedding $S_n \lra S_r \times \osym_c$ introduced in Corollary \ref{cor:real_toremb}, with $\mdimr G = r + 2c$. By restriction to $\Z_{m, 0}^n \cong \Z_m^{n - 1}$, Theorem \ref{thm:min_abl} guarantees that $c \ge n - 1$. Moreover, if $r < 2$, the composition $S_n \lra S_r \times \osym_c \lra \osym_c$ is an embedding and $c \ge n$ by the above. It follows that the minimal dimension is at least $2n$, and constructing $2n$-dimensional representations is straightforward.
\\ \\
For $m > 2$ and $n = 2$, $\mdimr G_0 (n, m) = 2$, as $G_0(2, m)$ is isomorphic to the dihedral group $D_{2m}$ on $m$ vertices. The dihedral groups provide examples where $\mdimc D_{2m} = \mdimr D_{2m}$ as they admit minimal faithful complex representations whose associated generating character subset $\fx$ is also conjugation-invariant. In other words, their minimal faithful representations are in fact minimal faithful real representations. 
\\ \\
The case $n = 4$ is particularly interesting, and we present it in a separate example. We also remark that $n = 2$ and $n = 4$ are distinguished largely because $S_{2} \cong \osym_{1}$ and $S_{4} \cong \osym_{3}$.
\end{example}

This last example shows the relationship between minimal real and complex dimension is more subtle than Proposition \ref{prop:rubiks} might suggest.  The minimal real representation constructed in Proposition \ref{prop:rubiks} came from taking a minimal complex representation and adding the necessary conjugate pieces or characters.  One might suspect this behavior is generic and that every minimal real representation comes from \emph{real-ifying} some minimal complex representation.  We conclude this work with Example \ref{ex:weirdexample} which shows this is not always the case.  A statement relating $\mdimc G$ and $\mdimr G $ in a meaningful way remains mysterious, and it would be of interest to say something precise for a family of groups.  We remark this example comes from the case by case analysis in Example \ref{ex:generalized_rubiks} and exists in part because $S_{4}$ is the unique symmetric group with a proper normal 2-subgroup.


\begin{example}\label{ex:weirdexample}
For each $m > 2$, the groups $G_{0}(4,m) = \Z_{m,0}^{4}\rtimes S_{4}$ where $S_{4}$ acts by the standard permutation action have $\mdimc G_{0}(4,m) = 4$. We find $\mdimr G_{0}(4,m) = 6$ for odd $m$ and $8$ for even $m$. For odd $m$, no minimal faithful real representation of $G$ is the `realification' of any minimal faithful complex representation of $G$.  
\end{example}

\begin{proof}
The first part of the above statement regarding the minimal complex dimension was proven in Proposition \ref{ex:generalized_rubiks}.  We now prove the claim about the minimal real dimension for odd $m$. For each $1 \leq i \leq 4$ define $\chi_{i}: \Z_{m,0}^{4} \lra \Z_{m} \leq S^{1}$ as the projection onto the $i$-th coordinate.  Next consider the collection $\fx = \{\chi_{2}^{-1}\chi_{3}^{-1},\chi_{2}\chi_{3},\chi_{1}^{-1}\chi_{3}^{-1},\chi_{1}\chi_{3},\chi_{1}^{-1}\chi_{2}^{-1},\chi_{1}\chi_{2}\}$.  Using the fact that $\chi_{1}\chi_{2}\chi_{3}\chi_{4} = 1$, this collection may be seen to be both $H$-invariant and generating.  It is also conjugation invariant, as $\chi \ol{\chi} = 1$ for every character of a finite abelian group.  Minimality follows because $\Z_{m,0}^{4} \cong \Z_{m}^{3}$ and by Theorem \ref{thm:min_abl}, every generating set of the dual of $\Z_{m,0}^{4}$ has at least $(2\cdot 3)$-elements, thus $\mdimr G  = 6$.  \\
\\
The collection $\fx$ fails to generate when $m$ is even, as the element $(m / 2, m / 2, m / 2, m / 2) \in \Z_{m, 0}^{4}$ lies in the kernel of each character. It follows that $\mdimr G > 6$ in this case, since $\fx$ is the only $H$-invariant and conjugation invariant collection of six complex characters, up to an overall exponent $\chi \mapsto \chi^k$. Moreover, one cannot adjoin a single real character to obtain a seven-dimensional representation, since real characters are permuted separately and the identity $\chi_1 \chi_2 \chi_3 \chi_4 = 1$ ensures no characters are fixed by $H$. Eight-dimensional representations follow readily from realification.
 \\ \\
For any $m > 2$, every minimal faithful complex representation of $G$ has $4$ complex characters, so the realification of any minimal faithful complex representation of $G$ will be $8$-dimensional, not $6$-dimensional.  Explicit calculations for $m = 3$ can be found on the first author's webpage \cite{e27s4}.  
\end{proof}

To conclude, we discuss some implications for orbifold singularities of minimal dimension. By Theorem \ref{thm:rmdim}, the minimal orbifold singularity of groups $G = AH$ locally resembles the quotient of a subgroup of $\Z_d^n \rtimes S_n$ acting on $\R^n$. The action of $H$ on $\hat{A}$ determines the minimal dimension $n$. Moreover, Corollary \ref{cor:when_ac} determines precisely which of these groups have orbifold singularities that come from restricting a complex manifold to its real structure. By Example \ref{ex:generalized_rubiks}, the groups $G_0(n, m)$ with $m > 2$ and $n \ne 2, 4$ admit minimal singularities which arise this way, whereas $G_0(2, m)$ and $G_0(4, m)$ for odd $m$ do not. The first follows from the isomorphism $G_0(2, m) \cong D_{2m}$, since the minimal real representations of $D_{2m}$ are non-orientable and thus cannot occur from the restriction of a complex structure. The second follows from Example \ref{ex:weirdexample}. Remarkably, this example provides orbifold singularities of dimension 6 which are orientable, but which fail to arise from restricting a complex orbifold singularity to a real one. 

\bibliographystyle{amsplain}
\bibliography{MinimalArxiv}

\end{document}